\documentclass[opre,nonblindrev]{informs3}
\usepackage{subfigure}
\OneAndAHalfSpacedXI 

\usepackage[ruled,norelsize,linesnumbered]{algorithm2e}

\DeclareMathOperator{\EX}{\mathbb{E}}
\usepackage{setspace}

\usepackage{endnotes}
\let\footnote=\endnote
\let\enotesize=\normalsize
\def\notesname{Endnotes}%
\def\makeenmark{\hbox to1.275em{\theenmark.\enskip\hss}}
\def\enoteformat{\rightskip0pt\leftskip0pt\parindent=1.275em
  \leavevmode\llap{\makeenmark}}

\usepackage[natbibapa]{apacite}
\usepackage{natbib}

\TheoremsNumberedThrough     
\ECRepeatTheorems

\EquationsNumberedThrough    

\begin{document}




\TITLE{Stochastic Dynamic Pricing for Same-Day Delivery Routing}

\ARTICLEAUTHORS{%
\AUTHOR{Anatolii Prokhorchuk, Justin Dauwels}
\AFF{School of Electrical and Electronic Engineering, Nanyang Technological University, Singapore 639798, \EMAIL{anatolii001@e.ntu.edu.sg\thanks{Corresponding author}, jdauwels@ntu.edu.sg}, 
} 
\AUTHOR{Patrick Jaillet}
\AFF{Department of Electrical Engineering and Computer Science, Operations Research Center, Massachusetts Institute of Technology, Cambridge, MA 02139, \EMAIL{jaillet@mit.edu}}
} 

\ABSTRACT{%
Same-day delivery for e-commerce has become a popular service. Companies usually offer several time delivery options with the earliest one being next hour delivery. Due to tight delivery deadlines and thin margins, companies often find it challenging to provide efficient same-day delivery services. In this work, we propose a holistic scheme that combines the optimization of routing and pricing for same-day delivery. The proposed approach is able to take into account uncertainty in travel times, a crucial factor for delivery applications in urban environments. We model this problem as a Markov decision process. We apply a value function approximation technique to compute opportunity costs. Based on these opportunity costs, as well as the customer choice model and travel time distribution, we optimize the prices for various delivery deadlines. We perform extensive computational experiments to compare the proposed model with baseline policies. We also investigate how the (potentially wrong) choice of travel time distributions affect the performance of the proposed optimization scheme. Through numerical simulations of realistic scenarios, we observe that compared to the deterministic model, the proposed approach can reduce the number of missed deliveries up to 40\%; at the same time, it can increase revenue by more than 5\% compared to the baseline policies. We explore new issues that arise due to the stochastic nature of the problem such as the effect of penalties for missed deliveries on pricing structure and overall revenue.
}%

\KEYWORDS{Routing; Pricing; Dynamic programming; Same-day delivery}


\maketitle


%

\section{Introduction}
According to a recent report by Technavio \citep{SDDreport}, Same-Day Delivery (SDD) market is expected to exceed USD 987 million with compound annual growth of over 90\% by 2020. Currently, same-day delivery is being offered by an increasing number of e-commerce companies (Amazon, Alibaba, Instacart, etc). Some provide deliveries within a 4-hour period, however, with the introduction of Amazon Prime Now, the delivery deadline has now been reduced to 1 hour \citep{primenow}. Instacart partners with physical stores and provides 1- and 2-hour deliveries as well. Given these tight time spans and small profit margins, it is not surprising that retailers struggle with making profits, especially in the food delivery sector \citep{FTdelivery}.
\par 
Same-day delivery exhibits certain characteristics that differentiate it from common vehicle routing problems. Namely, customer requests arrive dynamically during the day, therefore the delivery has to be made within a short period of time from a fixed depot. Each vehicle has to make several tours from and to depot during the shift, while any future route can be updated multiple times due to new customers' arrivals. \cite{voccia2017same} define this as a Same-Day Delivery Problem (SDDP) with the objective of maximizing the expected number of requests that can be delivered on time. In this problem, the pricing decisions are not considered.
\par
\begin{figure}[h!]
    \centering
     \subfigure{\includegraphics[width=0.87\columnwidth]{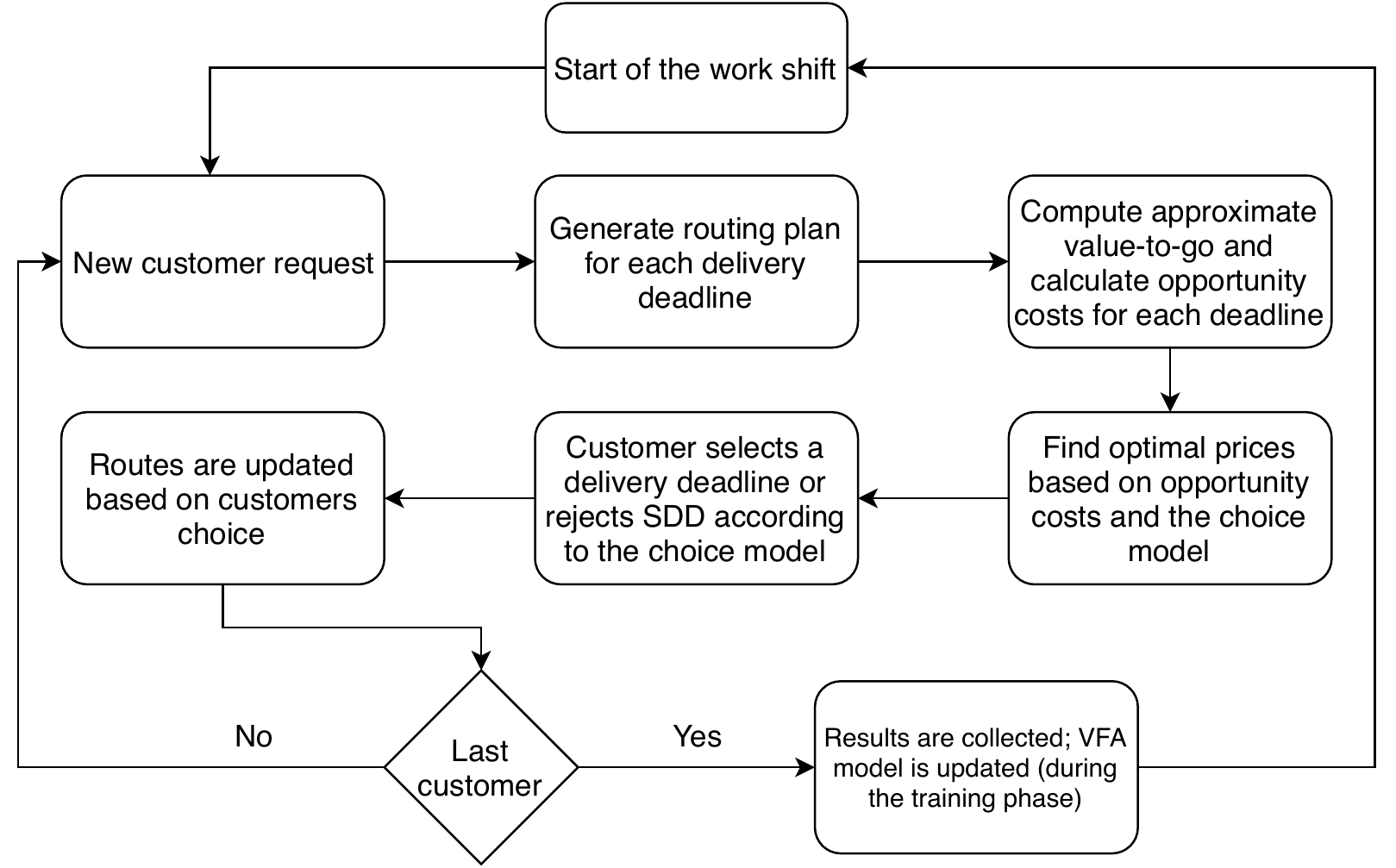}}
\label{figbs}
\caption{An overview of the model.}
\end{figure}

 
In real situations, delivery companies are also concerned with optimizing their profits, especially in such a competitive market. Pricing decisions affect both the current revenue and long-term customer loyalty. Modifying delivery prices to balance supply and demand is commonly employed in Attended Home Delivery (AHD) \citep{asdemir2009dynamic,yang2017approximate}. In this context, various methods have been proposed to provide different pricing for different time slots. \citet{ulmer2017dynamic} propose an approach for pricing same-day delivery options. The method computes opportunity costs via Value Function Approximation (VFA) for the Markov Decision Process (MDP). The existing literature on same-day delivery routing \citep{voccia2017same, ulmer2017dynamic, klapp2018dynamic}, however, only considers deterministic travel times and hence routing. We propose a method that takes into account the uncertainty of travel times. We similarly employ VFA to approximate the opportunity costs. However, the pricing decisions are based on solving a holistic optimization problem that considers the customer choice model, probabilities of arriving on time and possible penalties for late deliveries. The routing heuristic is also designed to maximize the expected revenue. The overview of the proposed offline approximate dynamic programming (ADP) algorithm is presented in Figure~\ref{figbs}.

\par 
Our contributions are as follows. We introduce a model called Stochastic Dynamic Pricing and Routing problem for Same-Day Delivery (SDPRSDD) where the stochasticity comes from both travel times and customer requests. We propose an approach based on value function approximation to solve the problem. To the best of our knowledge, this is the first study that deals with stochastic travel times in the context of same-day delivery routing and pricing. This problem involves a number of aspects not encountered in deterministic SDD routing. Namely, in addition to the commonly employed metrics such as the number of serviced customers and the overall revenue, we have to consider the number of customers for which the delivery arrived later than the deadline. We have to also consider the monetary effect of the missed deadline, specifically, whether the company compensates the delivery fee to the customer, or the company is additionally penalized for such deliveries (i.e. by giving a voucher to a customer). We perform an extensive set of computational experiments. We investigate the effect of incorporating the travel time distribution information on the overall performance. Since travel time information is often only known approximately, we also explore the effect of misspecifying the distribution on the performance. We compare different pricing policies and investigate the effect of different pairings of pricing policies and travel time information. The experiments cover a range of travel time distributions, customers spatial distribution, number of orders, the size of the fleet, as well as the penalties for missed deliveries.

\par 
The remainder of this paper is organized as follows. We present a literature review in Section 2 and formulate the problem definition in Section 3. We describe our proposed solution procedure in Section 4. We describe the computational experiments and the results in Section 5. We report our conclusions in Section 6.

\section{Relevant Literature}
There exists extensive literature on both dynamic delivery routing and pricing. However, only a few studies combine the two problems as mentioned by \cite{yang2014choice}. Here, we describe current approaches for joint routing and pricing optimization; further in the section, we provide a literature overview of the related studies on either routing or pricing. Table~\ref{tablelit} provides an overview of the most relevant approaches.
\subsection{Pricing and Routing} 
\cite{figliozzi2007pricing} describe a Vehicle Routing Problem in a Competitive Environment (VRPCE) where the provider bids on customer requests/contracts. The price is determined based on the expected loss in revenue (opportunity costs) which is computed via an online one-step look-ahead algorithm. The arrival time and the contract characteristics are not known in advance, however, the travel and service times are assumed to be deterministic. In contrast with the SDD problem, in VRPCE there is a need to provide only one price per order, and the acceptance of the contract is based on the auction. \cite{topaloglu2007incorporating} solve a truck dispatching problem by determining the prices beforehand. In this problem carrier prices influence transportation demand and the empty re-positioning costs are nonzero. 
The total expected profit is maximized by adjusting the variables following the sample-based directional derivatives of the objective function. The study that is most similar to ours is by \cite{ulmer2017dynamic}, where the author presents an approach for solving dynamic routing and pricing MDP for same-day delivery (DPPSDD). Value function approximation is applied to approximate the opportunity costs of accepting a customer. The pricing is then computed based on these opportunity costs (similarly to \cite{figliozzi2007pricing}). \cite{ulmer2017dynamic} employs meso-VFA which is a combination of parametric and non-parametric VFA described in \cite{ulmer2019meso}. Their numerical experiments show that the proposed method significantly outperforms baseline policies. 
\par 
In our study, we relax the assumption of deterministic travel times to make the model more applicable to real-world scenarios, however, this brings a new set of challenges. The first one is the fact that in a stochastic context it is impossible to guarantee the service within the deadline for most cases. In other words, the probability of arrival on time to serve a customer  cannot be guaranteed to be 1. Therefore, we also have to model the situations, in which the customer has not received the delivery before the guaranteed deadline. We consider a logit choice model for modeling customer behavior and incorporate it into the pricing optimization explicitly.
\par
Other research on delivery routing and pricing is mostly related to attended home delivery problem. Pricing policies are also a key aspect in ridesharing networks. We describe both applications further in this section. 
\par Next, we provide a short overview of relevant studies. We first describe previous research on stochastic routing with the focus on same-day delivery without the pricing considerations, then we mention other areas (AHD and ridesharing) where the pricing decisions are crucial.
\subsection{Stochastic and SDD Routing}
Stochastic Dynamic Routing \citep{bertsimas1991stochastic} is a well-studied topic. A detailed review of stochastic routing problems and approaches can be found, for example, in \cite{gendreau1996stochastic,adulyasak2015models,psaraftis2016dynamic,ritzinger2016survey,ulmerMDP}. However, only a limited number of studies deal with methods applicable to same-day delivery. This type of problem is characterized by short delivery time windows, vehicles performing multiple routes during a working shift and a necessity for a fast routing algorithm due to a dynamic request arrival. Related problems from a stochastic routing point of view include Vehicle Routing Problem with Stochastic Demands (VRPSD) \citep{bertsimas1990priori,bertsimas1996new}, Vehicle Routing Problem with Stochastic Customers (VRPSC) \cite{gendreau1995exact}, stochastic vehicle routing problem with deadlines (SVRP-D), \cite{adulyasak2015models}.  From the dynamic routing angle, previous work includes \citep{powell1996stochastic, regan1996dynamic} which deal with minimizing costs for trucking companies by deciding whether to accept or reject customer request.
\par 
There are several examples of more recent approaches that are similar to SDD routing. \cite{ghiani2009anticipatory} consider a vehicle dispatching problem with pickups and deliveries (VDPPD). They propose a sampling approach in which the number of samples is determined via Indifference Zone Selection (IZS) algorithm. This sampling procedure is applied to determine the best solution. \cite{bent2004scenario} propose a multi-scenario approach for dynamic vehicle routing problem with time windows (VRPTW), in which routing variants are generated for various scenarios incorporating known and unknown requests. \cite{azi2012dynamic} propose an adaptive large neighborhood search with local search heuristic. For each new request, the model considers multiple possible scenarios to decide whether to accept it. \cite{voccia2017same} propose to solve deterministic multi-trip team orienteering problem with time windows (MTTOPTW) during each decision epoch. Similarly to \cite{azi2012dynamic} they employ multiple-scenario approach based on \cite{bent2004scenario}. \cite{sungur2010model} consider a Courier Delivery Problem (CDP) with uncertain service times. They propose a scenario-based approach, combining an insertion heuristic and the tabu search algorithm. However, since these studies consider problems that differ from SDD, they usually do not consider such aspects as the customer choice model and pricing decisions.
\par 
Some of the most recent literature explicitly solves the SDD routing problem. \cite{klapp2016one, klapp2018dynamic} formulate Dynamic Dispatch Waves Problem (DDWP) for same-day delivery. In this problem, the objective is to determine whether to dispatch a vehicle during each decision step (`wave') or let it wait at the depot. For the deterministic version of the problem, they propose an integer programming approach. For cases when the customer arrival times are not known in advance, they propose several heuristics based on the deterministic solution. \cite{ulmer2018offline} propose a combination of VFA and online rollout algorithms to improve routing performance for single-vehicle routing problem with stochastic service requests.  \cite{ulmer2018preemptive} propose an approach based on approximate dynamic programming for SDD routing which allows preemptive returns to the depot.  \cite{yao2019} aim to solve a robust optimization SDD under the demand uncertainty. They model the problem as a precedence-constrained asymmetric TSP. To handle the computational complexity issues a mixed integer optimization approach is proposed. This approach outperforms the deterministic baseline. \cite{ulmer2018value} investigate how postponing customer requests from SDD to the next day affects the performance. To model these situations, they introduce a dynamic multi-period vehicle routing problem with stochastic service requests. By employing VFA with state space segregation and period classification they are able to outperform previous policies. Another application of VFA is described by \cite{van2017delivery}. Here, value function approximation is employed to solve the delivery dispatching problem with time windows. A simpler version of the problem is considered, in which the delivery time windows are replaced with the dispatch time windows. In addition to VFA, to allow the model to solve larger instances, an integer linear program is formulated to use within the ADP. \cite{ulmer2018same} study how drones can be combined with ordinary delivery vehicles for SDD. They propose a policy function approximation based on geographical districting to determine whether a particular order should be delivered with a drone or with a car. Finally, \cite{ulmer2019same} consider an SDD problem with pickup stations and autonomous vehicles. In this problem, goods are first consolidated and delivered to the pickup stations from the depot and then a fleet of autonomous vehicles performs the last-mile delivery. The problem is modeled as a Markov decision process and solved via a policy function approximation approach.
\par 
However, in all of the current SDD routing literature, the travel times are assumed to be deterministic in contrast with our approach. 

\begin{table}[]\scriptsize
\begin{tabular}{p{20mm}p{20mm}p{20mm}p{25mm}p{20mm}p{20mm}p{20mm}}
\hline
                            & \begin{tabular}[c]{@{}l@{}}Dynamic \\  Routing\end{tabular}                          & Travel Times & Pricing                            & Choice Model                                 & \begin{tabular}[c]{@{}l@{}}Time \\  Constraint \end{tabular}            & Model         \\                                                               \hline
Study                       &                                           &               &                                    &                                                    &                            &                                                                             \\ \hline
\cite{ulmer2017dynamic}                & Cheapest insertion                        & Deterministic & Opp. costs                  & Maximizing the utility based on WTP & SDD deadlines              & MDP + VFA                                                                   \\ \hline
\cite{figliozzi2007pricing}     & -                                         & Deterministic & Opp. costs                  & Auction                                            & Time windows               & \begin{tabular}[c]{@{}l@{}}Simulation (one-\\ step-look-ahead)\end{tabular} \\ \hline
\cite{topaloglu2007incorporating} & -                                         & Deterministic & Path-based directional derivatives & Demand depends on the price                        & -                          & Policy search                                                               \\ \hline
\cite{voccia2017same}        & Neighborhood search & Deterministic & -                                  & -                                                  & Time windows and deadlines & MDP + sampling                                                              \\ \hline
\cite{yao2019}          & MIO                                       & Deterministic & -                                  & -                                                  & -                          & MIO + uncertainty sets                                                      \\ \hline
\cite{Klein2018}        & -                                         & Deterministic & Opp. Costs                  & MNL                                                & Time slots                 & MILP                                                                        \\ \hline
\cite{yang2017approximate}     & -                                         & Deterministic & Opp. Costs                  & MNL                                                & Time slots                 & Insertion heuristic + historical demand data                                \\ \hline
\cite{ulmer2018offline}          & VFA + rollout algorithms                  & Deterministic & -                                  & -                                                  & -                          & MDP + VFA + rollout algorithms                                              \\ \hline
This study                  & Cheapest insertion                        & Stochastic    & Opp. Costs + customer model & MNL                                       & SDD deadlines              & MDP + VFA                                                                   \\ \hline
\end{tabular}
\caption{Classification of relevant studies.}
\label{tablelit}
\end{table}

\subsection{Attended Home Delivery}
Pricing and routing are commonly optimized together in the context of the attended home delivery. In AHD the goal is to provide pricing for time slots up to several days in advance. In that case, the final routing for a given day is usually fixed before the day starts which differentiates AHD from the SDD problem. Compared to dynamic routing, AHD optimization commonly concentrates on minimizing the routing costs instead of maximizing the expected number of served customers. \cite{campbell2005decision} propose a method to accept or reject a customer based on opportunity costs with routing approximated via insertion heuristic. In later work, \cite{campbell2006incentive} describe a model which introduces incentives for customers for choosing a specific time slot. \cite{koch2017} propose an approximate dynamic programming (ADP) approach for integrated pricing and routing problem. They introduce `time window budgets' concept and employ the least squares approximate policy algorithm to perform VFA. \cite{yang2014choice} propose approximating opportunity costs based on the insertion heuristic while also incorporating information about predicted demand. However, these predictions are based only on the historical data, without the information about the currently accepted customers. In later work, \cite{yang2017approximate} address some of the drawbacks by employing approximate dynamic programming. \cite{Klein2018} further improve upon this work by proposing a novel mixed-integer linear programming (MILP) approach that incorporates anticipation of the future demand, while the pricing is still computed via opportunity costs approximation.
\par 
Attended home delivery literature exhibits several similarities to SDD routing research: there is a need to provide pricing for delivery options and to consider customer choice model. However, since the routing can be planned up to several days in advance, the stochasticity of travel times is usually not considered at the point of customer request as the distributions are dependent on the current conditions.
\subsection{Ridesharing}
Properties of routing problems in ridesharing differ from the same-day delivery routing problem we consider. However, pricing decisions are crucial in both situations. \cite{furuhata2013ridesharing} provide an overview of the ridesharing approaches. In the context of this study, we are interested in dynamic pricing methods for ridesharing. One of the most known examples is the surge pricing employed by Uber; \cite{chen2016dynamic} discover that such an approach significantly increases the efficiency of the platform. \cite{spatpricing} investigate spatial price discrimination for a ridesharing platform. They observe that when the customer demand is not balanced, the optimal decision for the platform is to price rides differently according to the location. Applying this knowledge, we investigate how various pricing policies affect the service levels based on the location. We also look into how the demand spatial distribution (i.e. customer requests locations) influences the pricing discrimination. We investigate which policies result in `fair' decisions, e.g. whether policies always produce significantly higher prices for customers closer to the edge of the service area.

\section{Problem statement}

\begin{table}[htb]
\centering
\begin{tabular}{ll}
\hline
\multicolumn{1}{|l|}{Description}                  & \multicolumn{1}{l|}{Notation}                     \\ \hline
Work Shift                                         & $T_{\text{shift}} = [0 \ldots T] $                  \\
Last Order Time                                    & $t_{\text{last}}$                                 \\
Set of Vehicles                                    & $ \mathcal{V} = \{v_1, \ldots v_m \}$               \\
Set of Customers                                   & $ \mathcal{C} = \{C_1, \ldots C_K\}$               \\
Routing Plan                                       & $ \Theta = \{\theta(v_1), \ldots \theta(v_m)\} $ \\
Set of Delivery Deadlines                          & $ \Delta = \{\delta_1, \ldots \delta_d\}$         \\
Average delivery price for option $j$              & $P(\delta^j)$                                     \\
Probability of serving a customer within a deadline & $\mathbb{P}_R$                                    \\
Cost for violating the deadline (penalty)                & $c_{\text{miss}}$          
 \\
State of the MDP                                            & $S_k$   \\
Revenue                                            & $R$                                              
\end{tabular}
\caption{Notation.}
\end{table}

\subsection{Notation}
We employ the notation consistent with the DPPSDD problem from \cite{ulmer2017dynamic}. During the work shift $T_{\text{shift}} = [0, T]$, a fleet $\mathcal{V}$ of $m$ vehicles serves a set of dynamically arriving customers $ \mathcal{C} = \{C_1, \ldots C_K\}$. Customer arrival times and locations follow certain probability distributions. Vehicles make deliveries starting from depot $D$, while all customers are located within the service area $\mathcal{A}$. Each delivery order is first picked up at the depot. Given the short deadlines, and hence the short tours, we assume that the vehicles have unlimited capacity. Planned routing $\Theta = \{\theta(v_1), \ldots \theta(v_m)\}$ is a sequence of visits to customers and the depot for each vehicle: 
\begin{equation}
\theta(v_j) = ((C_{i_1}, a(C_{i_1}), \delta(C_{i_1})),(C_{i_2}, a(C_{i_2}), \delta(C_{i_2}))  \ldots (D, a(D), \delta(D)) \ldots  ),
\end{equation}
where $a(C_i)$ is the expected arrival time to the customer $C_i$ and $\delta(C_i)$ is the deadline. 
\par 
When the $k$-th customer $C_k$ request arrives at time point $t(C_k)$, the provider offers a set of delivery deadlines for that customer $\Delta^k = \{\delta_1, \ldots, \delta_{|\Delta^k|} \}$ with corresponding pricing vector $\mathcal{P}(C_k)$. $\delta_1$ corresponds to the next-day delivery which is always priced as $0$. Customer $C_k$ selects the preferred delivery deadline according to a logit choice model (see section 4.3 for details).

\subsection{Problem Statement}
We formulate the stochastic dynamic pricing and
routing problem for same-day delivery as a Markov decision process. We follow the route-based MDP formulation proposed in \cite{ulmerMDP} and investigated in \cite{ulmer2017dynamic}. The objective is to find a policy $\pi^{*}$ such as:
\begin{equation}
\pi^{*} = \argmax_{\pi} \left( \EX (\sum_{k=0}^K R(S_k,X^{\pi}(S_k))|S_0) \right) ,
\label{policy_eq}
\end{equation}
where $R$ is the revenue, $S_k$ is the $k$-th state and $X^\pi$ is the decision. The initial state $S_0$ occurs when $t=0$. Each customer request $k$ represents a decision point at time $t(C_k) := t_k$. In each decision point, the MDP state $S_k$ is defined as:  \begin{equation}
    S_k = (t_k, \Theta_k,\mathcal{C}_k,C_k).
\end{equation}
The state contains information about the current time, the currently planned routing, as well as information about all existing customers up to $k$-th (including the location, chosen deadline and price). Decisions in this formulation represent the pricing offer to the customer with the corresponding routing plans: $x_k = (\mathcal{P}_k, \mathcal{U}_k)$. Each element $P^j_k \in \mathcal{P}_k$ represents a price for a certain delivery option, $P^j_k \in \mathbb{R}_+$. There is also a corresponding routing plan for each option: $\mathcal{U}_k = \{ \Theta_k^j \}$. However, the decision does not necessarily contain prices for each delivery option. 
\par 
After making a decision, the next state is obtained via stochastic transition $\omega_k$. Similarly to \cite{ulmer2017dynamic}, the transition consists of the customer selecting one of the proposed options and vehicles following the proposed routing plans until the next customer request arrives. However, in our case, the travel times are not known in advance. The random instance of the travel time is only generated when a given vehicle is about to travel through the next arc on its path. However, to make the model more realistic, during the actual simulations the random values are generated once for each area and time period (the details are presented in Section 5).
\subsection{Example}
\begin{figure}[h!]
    \centering
     \subfigure{\includegraphics[width=0.47\columnwidth]{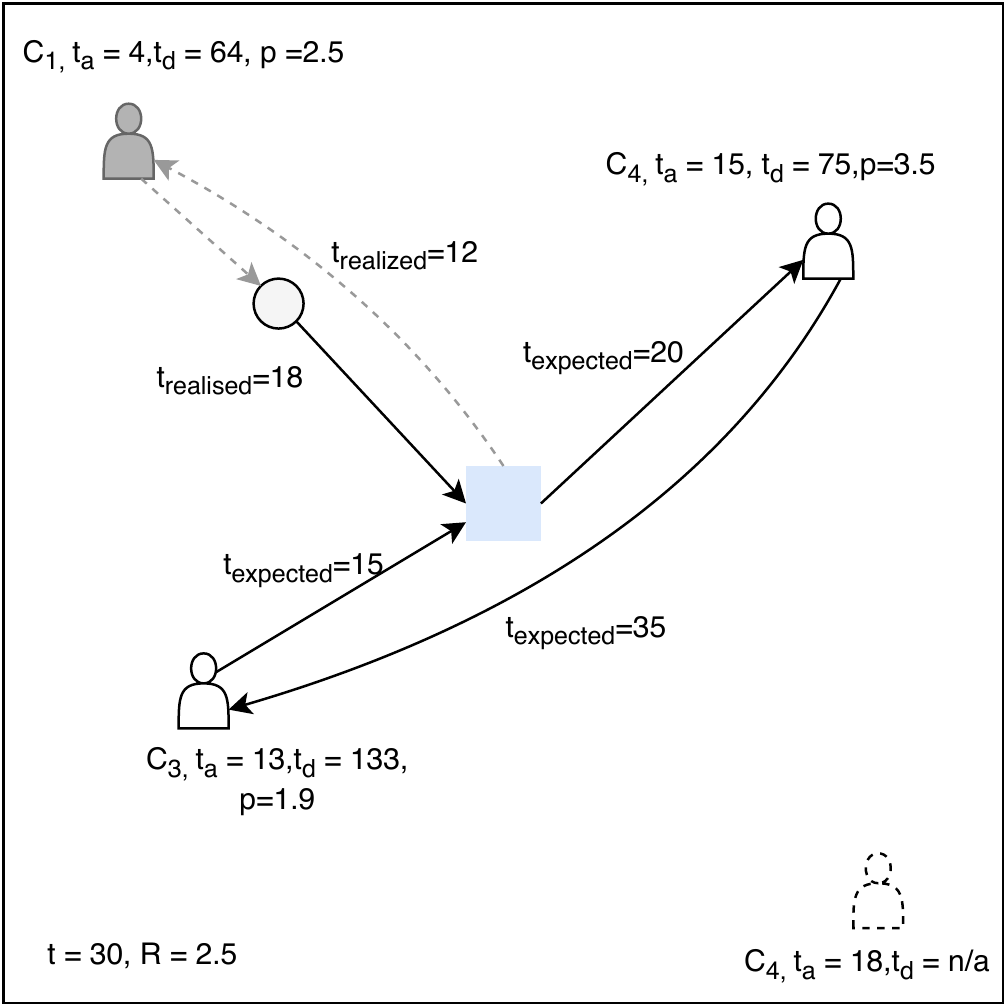}}
     \subfigure{\includegraphics[width=0.47\columnwidth]{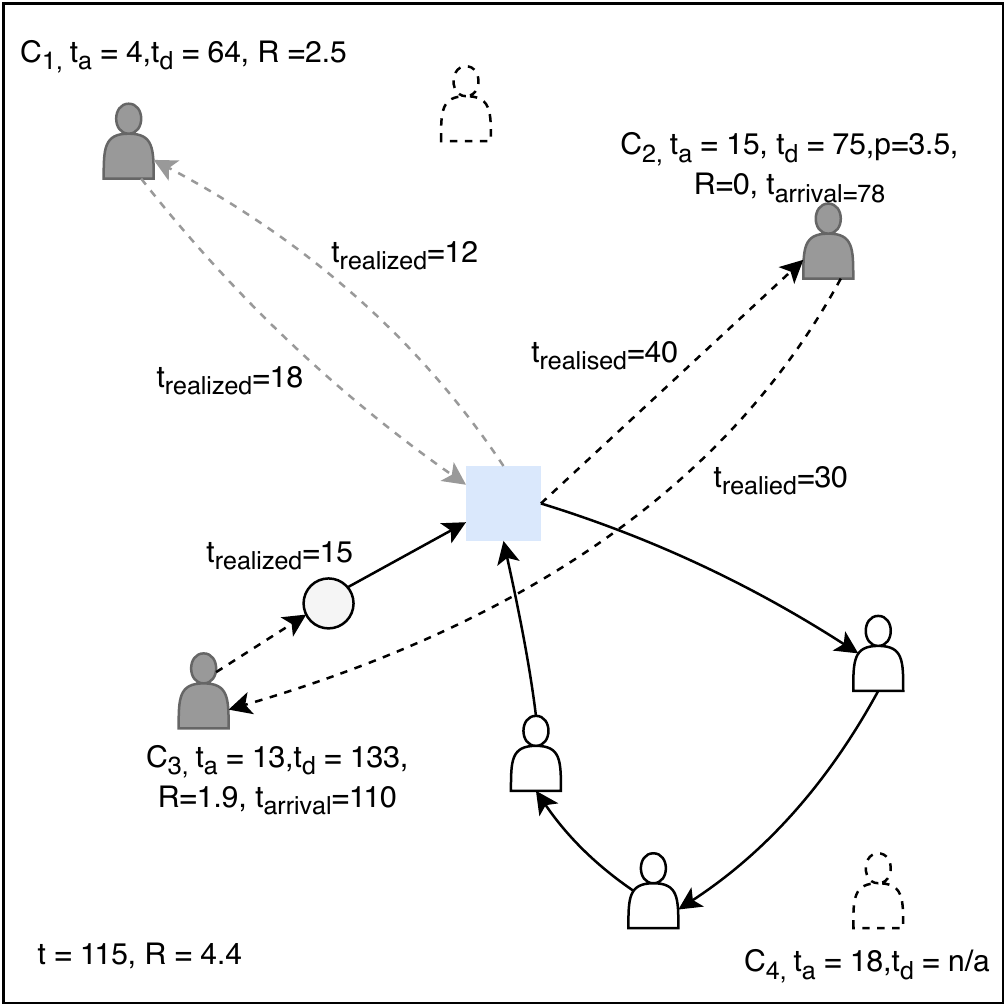}}
    \caption{Illustration of SDPRSDD at two time instances. The depot is represented as a square, finished routes as dotted arrows, planned routes as solid arrows, served customers are shaded, and customers who chose next-day delivery are dotted.}
    \label{example}
    \end{figure}

An illustration of the problem is shown in Figure~\ref{example}. For simplicity, we assume that only one vehicle services the depot. The diagram on the left corresponds to time $t=30$ which is between two decision points: $t(C_4)$ and $t(C_5)$. The requests from 4 customers have been received up to this point. The vehicle has already made a delivery to $C_1$ and the current planned route includes a stop in the depot and servicing customers $C_2$ and $C_3$. $C_4$ declined the same-day delivery and is therefore not in the route. The delivery price for $C_1$ is already added to the current revenue. The diagram on the right shows the same work shift at time $t=115$. Customers $C_2$ and $C_3$ have already been serviced, and 4 new customers have arrived. The travel time from the depot to $C_2$ was sampled to be equal $40$ which led to missing the delivery deadline to $C_2$. This delivery does not yield any revenue. The deadline for $C_3$ was satisfied, hence the overall revenue is increased by $1.9$.

\par 

\section{Approach}

\subsection{Bellman Equation}
To solve equation \ref{policy_eq} we employ the Bellman principle of optimality. At each state $S_k$ we can solve the equation by maximizing the sum of expected immediate revenue and the expected future rewards conditionally on $S_k$ and $x$:
\begin{equation}
    \argmax_{x} \EX R(S_k, x) + \EX \sum_{l=k+1}^{K}[R(S_l, X_l^{\pi^*}(S_l))|S_l,x ].
\end{equation} 
The immediate revenue can be calculated as follows: \begin{multline}
        \EX R(S_k, x) = \sum_{\delta^j \in \Delta} \mathbb{P}(\delta^j,S_k,\mathcal{P}) \Bigg[ \mathcal{P}(\delta^j) \times  \mathbb{P}_R(\delta^j, C_k, \hat{\Theta^j})  - c_{\text{miss}} \times (1 - \mathbb{P}_R(\delta^j, C_k, \hat{\Theta^j})) + \\ \sum_{ \substack{ C_i : i < k\\ C_i \text{ is not yet served}}} P(\delta(C_i)) \times \mathbb{P}_R(\delta(C_i),C_i, \hat{\Theta^j}) - c_{\text{miss}} \times (1 - \mathbb{P}_R(\delta(C_i),C_i, \hat{\Theta^j})) \Bigg],
\end{multline}
where $P(\delta(C_i))$ is a price the customer $C_i$ has agreed to pay for the delivery,  $\mathbb{P}_R(\delta(C_i),C_i, \hat{\Theta^j})$ is the probability of serving the customer $C_i$ within the deadline $\delta(C_i)$ based on the route $\hat{\Theta^j}$, and  $c_{\text{miss}}$ is a penalty that we are charged in case a particular customer is not served within agreed deadline. \par
The second part of the equation (4) is called value $V$:
\begin{equation}
    V(S_k,x_k) = \sum_{\delta^j \in \Delta} \mathbb{P}(\delta^j,S_k,\mathcal{P}) \times V(S_k, \hat{\Theta}^j),
\end{equation}
where the summation is over all possible deadlines $\delta^j$, and $\hat{\Theta}^j$ is a routing plan for each deadline.  
\subsection{Overview}
The proposed method is described in Algorithm 1. Here we provide a brief overview of the approach before describing each part in detail. First, for each new customer request and for each possible delivery deadline (lines \ref{alg:get_cust}:\ref{alg:deadline}), we compute the best routing plan via the cheapest insertion heuristic (line \ref{alg:routing}). For this problem, we consider three possible deadlines: 1, 2 and 4 hours. We select the best route based on the expected profit which can be calculated from the information about travel time distribution. If none of the insertion positions produces an increase in expected revenue for a particular delivery deadline, then it is considered infeasible and not offered to the customer. If all delivery deadlines are infeasible then the customer is only offered next-day delivery (which we assume they always accept in this situation). For each feasible routing plan, we compute the opportunity costs by value function approximation (lines \ref{alg:vfa2}:\ref{alg:opp}). Concretely, we perform VFA by linear regression of several features representing the current state. Next given the resulting opportunity costs and the customer choice model we find the optimal pricing vector by maximizing the expected revenue (line \ref{alg:opt_pricing}). The computed vector is then presented to a customer and the customer choice is sampled from the choice model (line \ref{alg:logit}). If the customer selected same-day delivery, we update the current routing plan, otherwise no changes are made.

\subsection{Routing}
The cheapest insertion heuristics is commonly \citep{ulmer2017dynamic, Klein2018} applied in problems combining routing and pricing due to its computational speed and reasonable performance. We propose to employ a modified version of this heuristic, where the cheapest insertions are based on the expected revenue of the route instead of the travel time extension. This allows us to apply this heuristic in a stochastic context. For a customer $C_k$, for each possible delivery deadline $\delta^j$ we select a route according to the following optimization problem:
\begin{multline}
    \argmax_{\Theta^j}  \Bigg[ P(\delta^i) \times  \mathbb{P}_R(\delta^j, C_k, \Theta^j)  - c_{\text{miss}} \times (1 - \mathbb{P}_R(\delta^j, C_k, \Theta^j)) + \\ \sum_{ \substack{ C_i : i < k\\ C_i \text{is not yet served}}} P(\delta(C_i)) \times \mathbb{P}_R(\delta(C_i),C_i, \Theta^j) - c_{\text{miss}} \times (1 - \mathbb{P}_R(\delta(C_i),C_i, \Theta^j)) \Bigg],
\label{exprev}
\end{multline}
where $P(\delta^j) := w^j$ is an average willingness to pay (or, an average price that is charged) across the population for a particular deadline $j$.
\par 
Solving (\ref{exprev}) gives us updated route $\hat{\Theta^j}$ and its expected revenue $\hat{r}_k^j$ for customer $C_k$ and deadline $\delta^j$. We are going to offer deadline $\delta^j$ for this customer only if $\hat{r}_k^j > \hat{r}_k$, where $\hat{r}_k$ is the expected revenue of the current routing (before the addition of customer $C_k$). In other words, we consider the delivery deadline to be feasible if the expected change in revenue is positive.

\subsection{Opportunity Costs}
Similar to \cite{ulmer2017dynamic} and \cite{Klein2018} we rely on approximations of opportunity costs to provide pricing options. 
\par 
From solving the routing problem in the previous section we obtain up to $3$ routing updates: $\{ \hat{\theta}^{60},\hat{\theta}^{120},\hat{\theta}^{240} \}$, which correspond to $60$-min, $120$-min and $240$-min deadlines. For each $\hat{\theta}$ we are going to compute the opportunity costs as:
\begin{equation}
    \hat{O}^j = \hat{V}(\theta^j) - \hat{V}(\theta^0),
\end{equation} 
where $\hat{V}(\cdot)$ is the value function and $\theta^0$ is the current route (before the update).
\par 
Due to the curse of dimensionality, it is intractable to compute the value function directly. Therefore, we utilize Value Function Approximation (VFA, \cite{powell2007approximate}) to compute $V(\theta)$. We employ parametric VFA and represent the value function as a linear function for each time period $\tau$:
\begin{equation}
    \hat{V}_\tau (\theta) = c_\tau^b + \sum_i f_i c_\tau^i,
\label{eq:vfa}
\end{equation}
where $c_\tau^b$ is a baseline coefficient, $f_i$ is a value of the $i$-th feature and $c_\tau^i$ is the corresponding coefficient.
\par 
The list of features include:
\begin{itemize}
    \item Free time budget: \begin{equation}
        f_1(\theta) = t_{\text{max}} - a(D_2^\theta),
    \end{equation} where $a(D_2^\theta)$ is the expected arrival time to the last stop in the route (the depot).
    \item Flexibility (proposed in \cite{ulmer2017dynamic}): 
    \begin{equation}
    f_2 =
    \left\{
    	\begin{array}{ll}
    		|\mathcal{C}_n(\theta)|^{-1} \sum_C (\delta(C)-a(C))  & \mbox{if } |\mathcal{C}_n(\theta)|\geq 0 \\
    		t_{\text{max}}-t_k & \mbox{else } 
    	\end{array},
    \right\}
    \end{equation}
    \item Probability that all customers will be served on time: \begin{equation}
        f_3(\theta) = \prod_i \mathbb{P}_R(C_i, \theta),
    \end{equation}
    \item Time budget in the `worst case'. Instead of relying on expected travel times, the budget is calculated with mean arrival time values plus 2 standard deviations: \begin{equation}
        f_4(\theta) = t_{\text{max}} - a(D_2^\theta)_{\text{2std}},
    \end{equation}
    \item Average distance per customer in the current delivery tour: 
    \begin{equation}
        f_5(\theta) = |\mathcal{C}_n(\theta)|^{-1} (\sum_{k=2}^n d(C_{k-1}, C_k) + d(D, C_1) + d(C_n, D)).
    \end{equation}
\end{itemize} 
\subsection{Pricing}
From sections 4.3 and 4.4, we obtain a set of up to 3 possible route updates and corresponding opportunity costs: $\{ (\theta^{60}, \hat{O}^{60}),(\theta^{120}, \hat{O}^{120}),(\theta^{240}, \hat{O}^{240}) \}$.

To obtain a pricing vector we solve the following optimization problem: 
\begin{equation}
    \mathcal{P}^{*}_k = \argmax_\mathcal{P} \sum_{\delta^j} \mathbb{P}_C(\mathcal{P}, \mathcal{P}(\delta^j)) \times \lbrack  (\mathcal{P}(\delta^j) \times \mathbb{P}_{\text{ontime}} - c_{\text{miss}} \times (1 - \mathbb{P}_{\text{ontime}}))   - \hat{O}^j \rbrack,
\end{equation}
where $\mathbb{P}_C(\mathcal{P}, \mathcal{P}(\delta^j))$ is the probability that customer will choose option $j$ given the pricing vector $\mathcal{P}$ and $\mathcal{P}(\delta^j)$ is the price corresponding to option $j$. This problem represents maximization of the expected revenue over all possible customer choices, taking into account the probability of the customer accepting a particular option $j$ and the probability of satisfying this deadline. The corresponding opportunity costs $\hat{O}^j$ are subtracted from each expected revenue term to account for the fact that satisfying any particular deadline means losing possible future value due to this capacity allocation. 
\par 
We employ the logit choice model to quantify customers behavior. In that model, the probability that the customer will select the option $j$ is given by:
\begin{equation}
    \mathbb{P}_C(\mathcal{P},\delta^j) = \frac{e^{u^k_j}}{e^{u^k_0} + \sum_{i} e^{u^k_i} },
\end{equation}
where $u^k_i$ is utility of customer $k$ and option $i$ and  $u^k_0$ is utility of next-day delivery, which we set to $0$. Since utilities of customers are unknown exactly due to stochasticity, during optimization we rely on the expected values to approximate probabilities. The deterministic optimization problem is then solved by the L-BFGS algorithm with bounds \citep{byrd1995limited}. In cases when one or more delivery options are not feasible, we set the price as  a large constant (10000), which results in zero probability of selecting these options.

\begin{algorithm}[H]
\caption{SDPRSSD}
\SetAlgoLined
 \For{$t\gets1$ \KwTo $T$}  { 
  \If{\text{ExistsRequest}(t)}{
   $C_k \gets$ GetCustomer(t)\; \label{alg:get_cust}
   $\Theta \gets$ CurrentRouting(t)\;
   $O \gets \varnothing$\;
   $\mathbb{P}_{\text{ontime}} \gets \varnothing$\;
   \ForAll{$\delta_j \in \Delta_{\text{SDD}}$}  { \label{alg:deadline} 
   $\hat{\Theta}^j \gets \text{InsertionHeuristicRouting}(\Theta,t, C_k, \delta_j)$ \; \label{alg:routing}
   $V_0 \gets \text{VFA}(t,\Theta)$ \;  \label{alg:vfa}
   $V_j \gets \text{VFA}(t, \hat{\Theta}^j)$ \; \label{alg:vfa2}
   $O^j \gets V_0 - V_j$ \; \label{alg:opp}
   $O \gets O \cup \{ O^j \}$ \;
   $\mathbb{P}_{\text{ontime}}^j \gets \text{ProbabilityOfArrivingOnTime}(\hat{\Theta}^j, \delta_j) $\; \label{alg:poba}
   $\mathbb{P}_{\text{ontime}} \gets \mathbb{P}_{\text{ontime}} \cup \{ \mathbb{P}_{\text{ontime}}^j \} $ \;
   }
   $\mathcal{P} \gets \text{OptimizePricing}(O, \mathbb{P}_{\text{ontime}})$  \; \label{alg:opt_pricing}
   $c \gets \text{LogitChoiceModel}(C_k, \mathcal{P})$ \; \label{alg:logit}
   \If{$c>1$} {
    $\Theta \gets \hat{\Theta}^c$ \;
   }
   }
   $\Theta \gets \text{UpdateRouting}(\Theta, t)$ \;

 }
\end{algorithm}

\section{Computational Experiments}
\subsection{Experiment Description}
We simulate a scenario similar to \cite{ulmer2017dynamic}. We assume the time limit $t_{\text{max}} = 480$, the last order time $t_{\text{last order}} = 420$, the loading and service times $t_{s}=t_{D}=2$ and possible delivery deadlines of $60$, $120$, and $240$ minutes. Customers arrive according to a Poisson process. The expected number of orders is different across different instances.
\par
We test two travel time distributions. Specifically, we assume that the inverse of the speed is distributed according to one of the following distributions:
\begin{itemize}
    \item Gaussian mixture with two components with equal weights: $\mathcal{N}(1/20,1/250)$ and $\mathcal{N}(1/40,1/250)$,
    \item Gaussian distribution $\mathcal{N}(0.0375,0.0131)$.
\end{itemize}
Both distributions have the same mean and variance. In both cases the sampled values are bounded between $\frac{1}{5}$ and $\frac{1}{120}$.

\par 
Distance between any two points is assumed to be Euclidean; we set km/h as a unit for speed and km as unit for distance. Generating speed values every time a vehicle starts the trip would be unrealistic: this can lead to situations when two vehicles depart the same location one after another and experience different travel times. To closer emulate real-world conditions, we instead generate one sampled speed value for all trips starting at particular area during particular time period. To this end, we divide the service area into 4 quadrants around the center (0,0) and set the time period to 15 minutes.
\par 
We compare the performance on three different spatial distributions for customers locations. In the first scenario, there are two independent Gaussian distributions: both $x$ and $y$ coordinates are distributed as $\mathcal{N}(0,2.5)$. The second model consists of uniform distributions for both coordinates: $x, y \sim  U(-10,10)$. In the third scenario, the customers are centered around 4 clusters. Cluster centers and corresponding means of Gaussian r.v. are located at $(5,5), (5,-5), (-5,5), (-5,-5)$, with standard deviation equal to $1$. In all cases, the depot is located in the center $(0,0)$.
\subsection{Customer Choice Model}
\begin{figure}[h!]
    \centering
     \subfigure{\includegraphics[width=0.6\columnwidth]{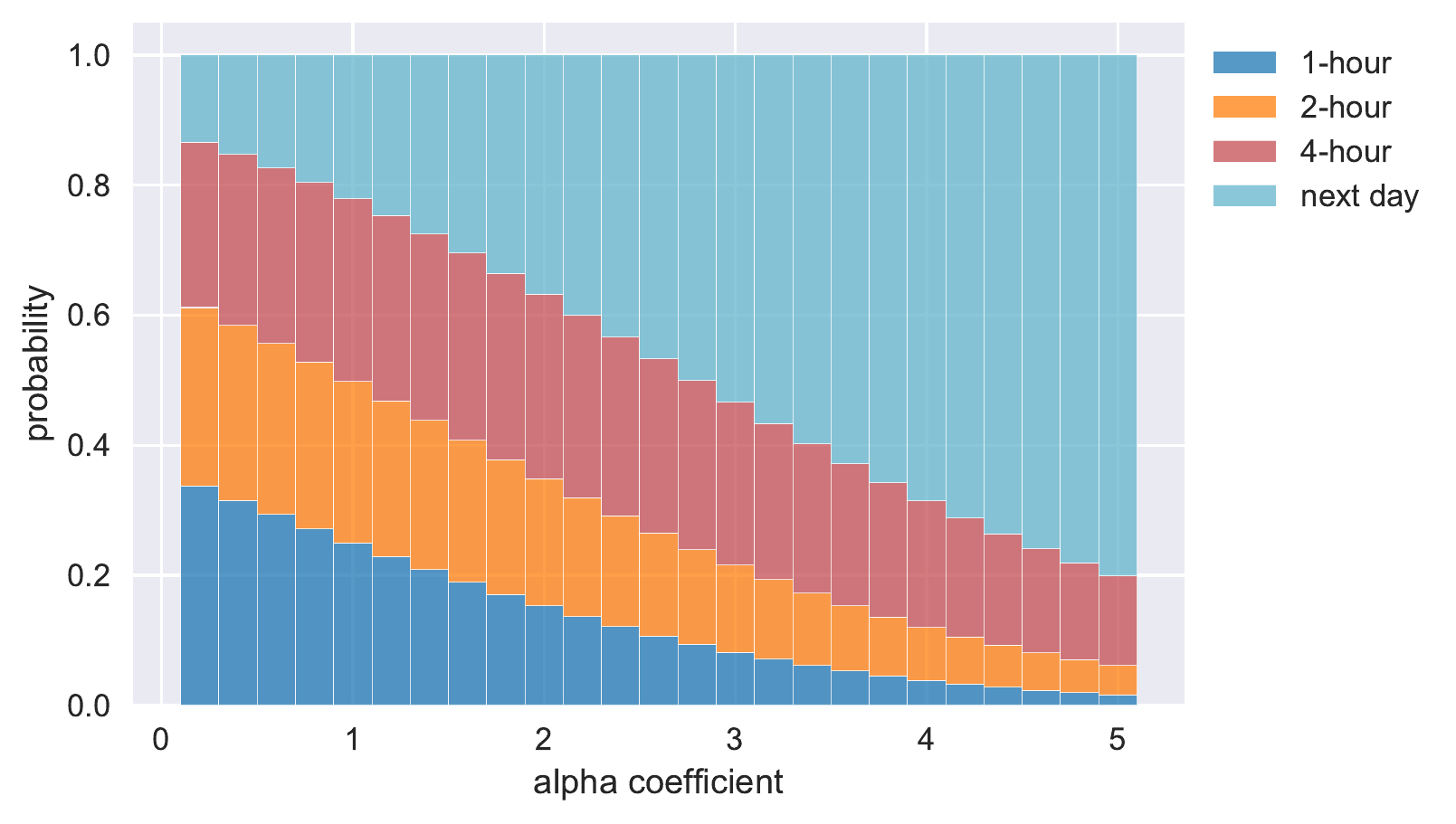}}
    \caption{Probability of accepting a particular delivery deadline according to the logit choice model based on pricing vector $[\alpha, 0.75\alpha, 0.5\alpha ]$ for different values of $\alpha$.}
    \label{figchoice}
    \end{figure}

The Multinomial Logit model (MNL) is a discrete-choice model in which the decision-maker chooses the option that maximizes their utility \citep{ben1985discrete}. This model is commonly applied in delivery problems, e.g. in \cite{yang2014choice, Klein2018}. We employ the logit choice model to simulate customer behavior. For each customer $C_i$ and delivery deadline $\delta_j$  utility is defined as:
\begin{equation}
    u_i^j = \beta_1 \delta^j + \beta_2 (\delta^j)^2 - p_j + \varepsilon_i^j,
\end{equation}
where $\beta_1, \beta_2$ are model parameters, $p_j$ is an offered price for the delivery option and $\varepsilon_i^j$ is Gumbel distributed noise.  We choose $\beta_l$ such that $\sum_l \beta_l (\delta^{60})^l = 1$, $\sum_l \beta_l  (\delta^{120})^l = 0.75$, $\sum_l \beta_l (\delta^{240})^l = 0.5$. Probabilities of choosing different pricing levels $\alpha$ with corresponding pricing vector $[\alpha, 0.75\alpha, 0.5\alpha ]$ are illustrated in Figure~\ref{figchoice}.

\begin{figure}[h!]
    \centering
     \subfigure{\includegraphics[width=0.6\columnwidth]{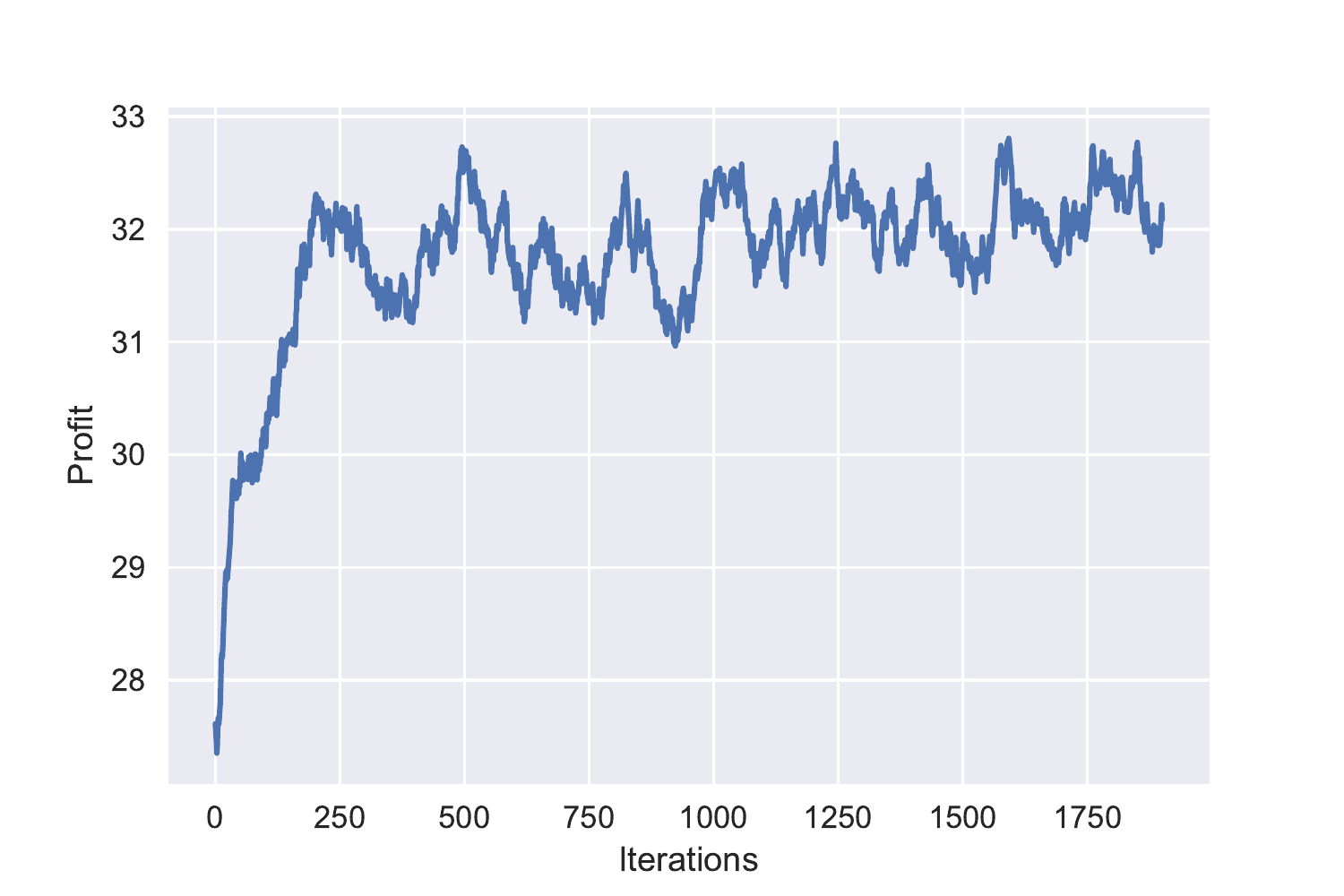}}
    \caption{Running average of profit during the training for one of the instances.}
    \label{figiter}
    \end{figure}

\begin{table}[]
\renewcommand{\arraystretch}{1.1}
\centering
\begin{tabular}{|l|l|l|l|l|l|l|}
\hline
Instance                                                                                             & 1-hour    & 2-hour    & 4-hour    & Next-day & Rejected & Revenue \\ \hline
\begin{tabular}[c]{@{}l@{}}80 orders, 1 vehicle, 0 penalty, \\ Gaussian customers\end{tabular} & 12.7\%/2.20 & 12.9\%/2.01 & 10.1\%/1.97 & 45.1\%/0 & 19.2\%/0 & 56.65 \\ \hline
\begin{tabular}[c]{@{}l@{}}80 orders, 1 vehicle, 2 penalty, \\ Gaussian customers\end{tabular} & 12.1\%/2.24 & 13.2\%/1.96 & 10.0\%/1.98 & 45.3\%/0 & 19.4\%/0 & 55.52 \\ \hline
\begin{tabular}[c]{@{}l@{}}80 orders, 1 vehicle, 0 penalty, \\ Uniform customers\end{tabular} & 5.2\%/3.00 & 8.8\%/2.52 & 11.1\%/2.16 & 62.1\%/0 & 12.8\%/0 & 37.19 \\ \hline
\begin{tabular}[c]{@{}l@{}}80 orders, 1 vehicle, 2 penalty, \\ Uniform customers\end{tabular} & 3.1\%/3.54 & 6.9\%/2.87 & 10.5\%/2.24 & 65.0\%/0 & 14.5\%/0 & 35.15 \\ \hline
\begin{tabular}[c]{@{}l@{}}80 orders, 3 vehicles, 0 penalty, \\ Gaussian customers\end{tabular}  & 20.7\%/1.94  & 16.7\%/1.94  & 12.8\%/1.94 & 49.7\%/0  & 0.1\%/0 & 77.26  \\ \hline
\begin{tabular}[c]{@{}l@{}}80 orders, 3 vehicles, 2 penalty, \\ Gaussian customers\end{tabular}  & 20.6\%/1.96  & 16.6\%/1.95  & 13.0\%/1.94 & 49.7\%/0  & 0.1\%/0 & 74.87   \\ \hline
\begin{tabular}[c]{@{}l@{}}80 orders, 3 vehicles, 0 penalty, \\ Uniform customers\end{tabular}  & 12.9\%/2.34  & 16.2\%/2.05  & 15.0\%/1.88 & 55.8\%/0  & 0.1\%/0 & 63.47   \\ \hline
\begin{tabular}[c]{@{}l@{}}80 orders, 3 vehicles, 2 penalty, \\ Uniform customers\end{tabular}  & 10.6\%/2.76  & 13.6\%/2.29  & 14.9\%/1.98 & 60.8\%/0  & 0.1\%/0 & 62.22  \\ \hline
\end{tabular}
\caption{Effect of different instance specifications on customer service levels. Values in the cells are average proportion of customers selecting a particular deadline / average price the customer paid for a particular deadline. In all cases, the travel time distribution and the model distribution are Gaussian.}
\label{simple_comp}
\end{table}

\begin{figure}[]
    \centering
     \subfigure{\includegraphics[width=0.47\columnwidth]{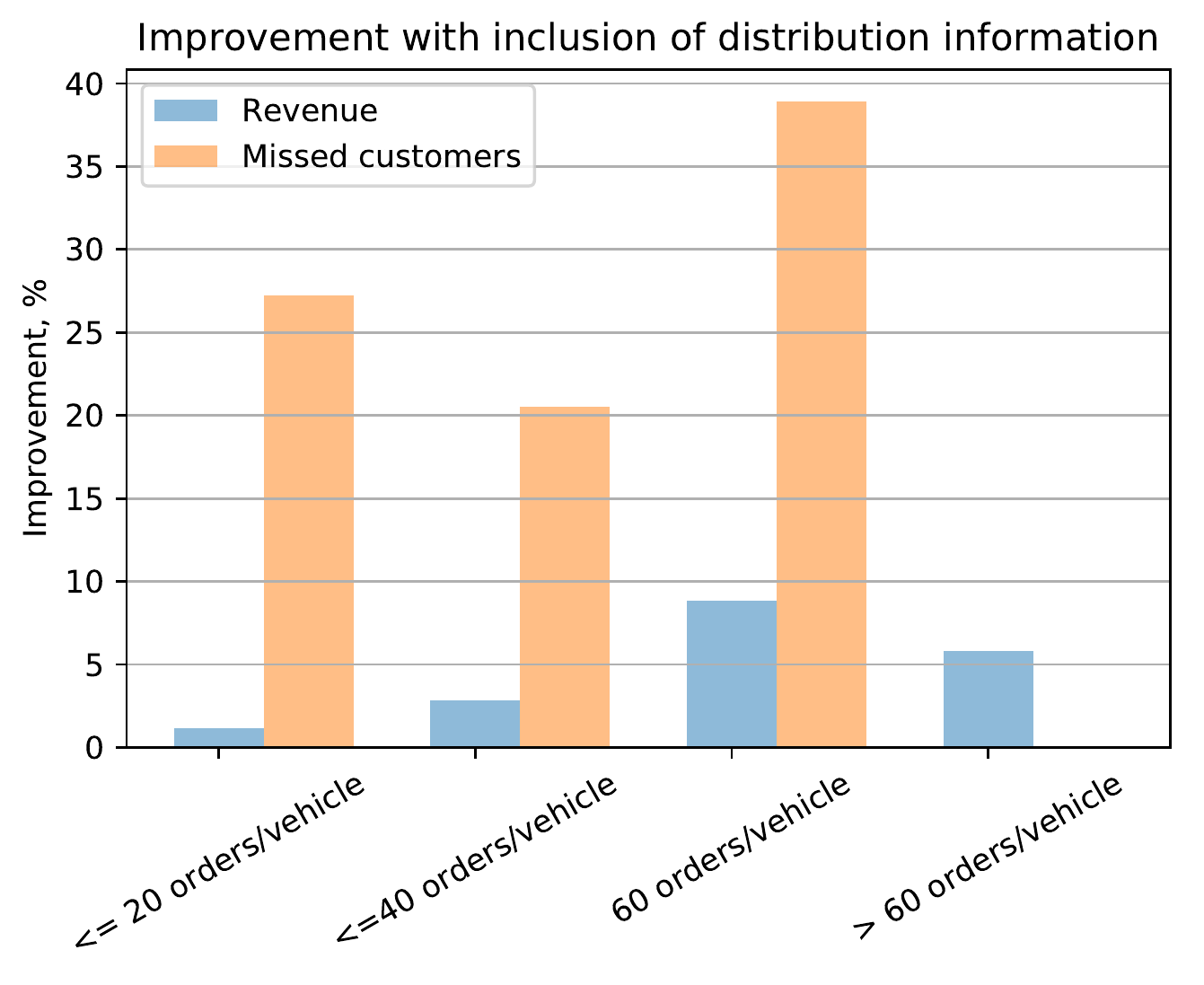}}
     \subfigure{\raisebox{13pt}{\includegraphics[width=0.47\columnwidth]{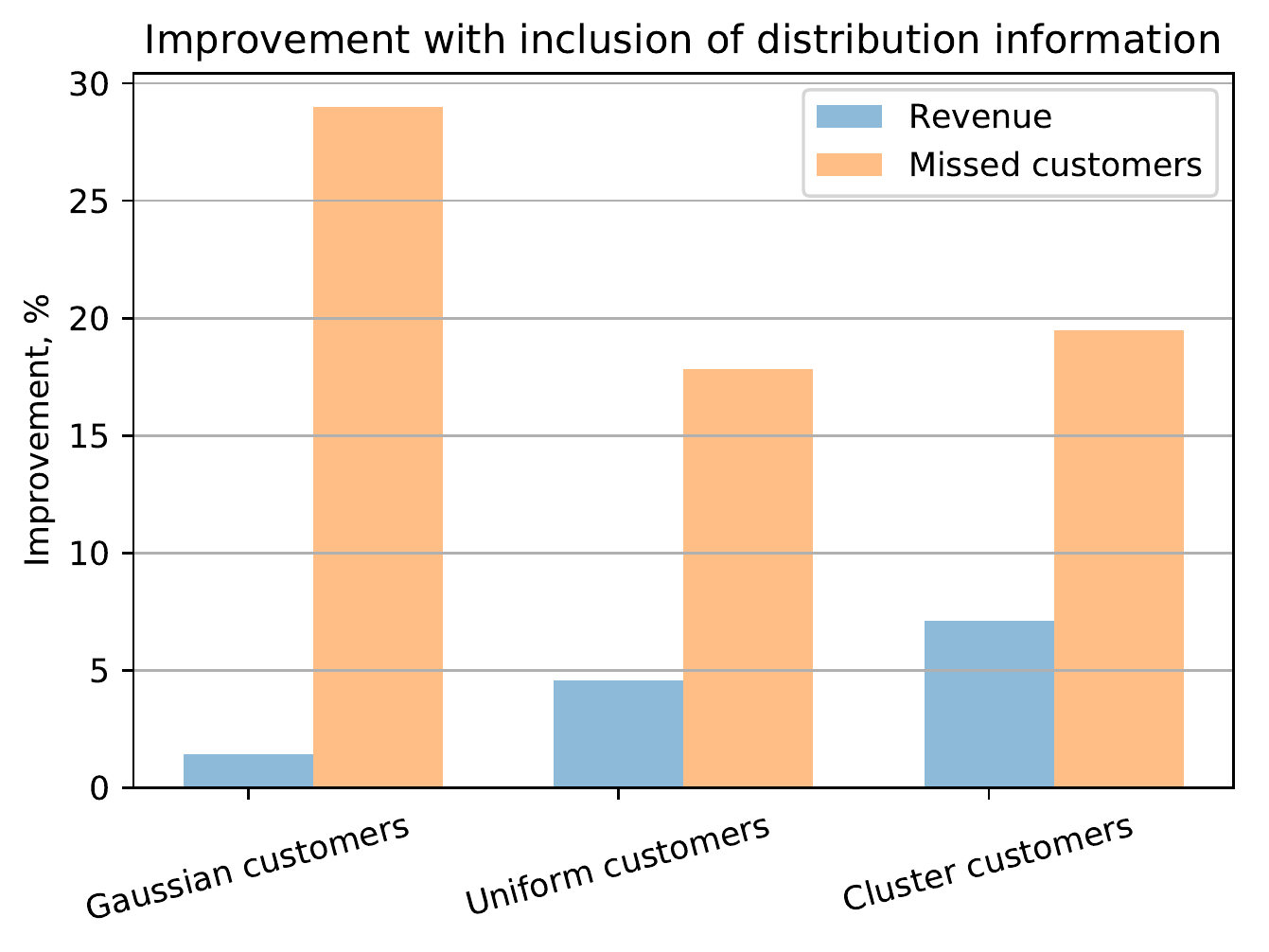}}}
    \caption{Percentage improvement in expected revenue and decrease in missed customers due to incorporating the travel time distribution information across different segments.}
    \label{res_imp}
    \end{figure}
    
    \begin{figure}[]
    \centering
     \subfigure{\includegraphics[width=0.47\columnwidth]{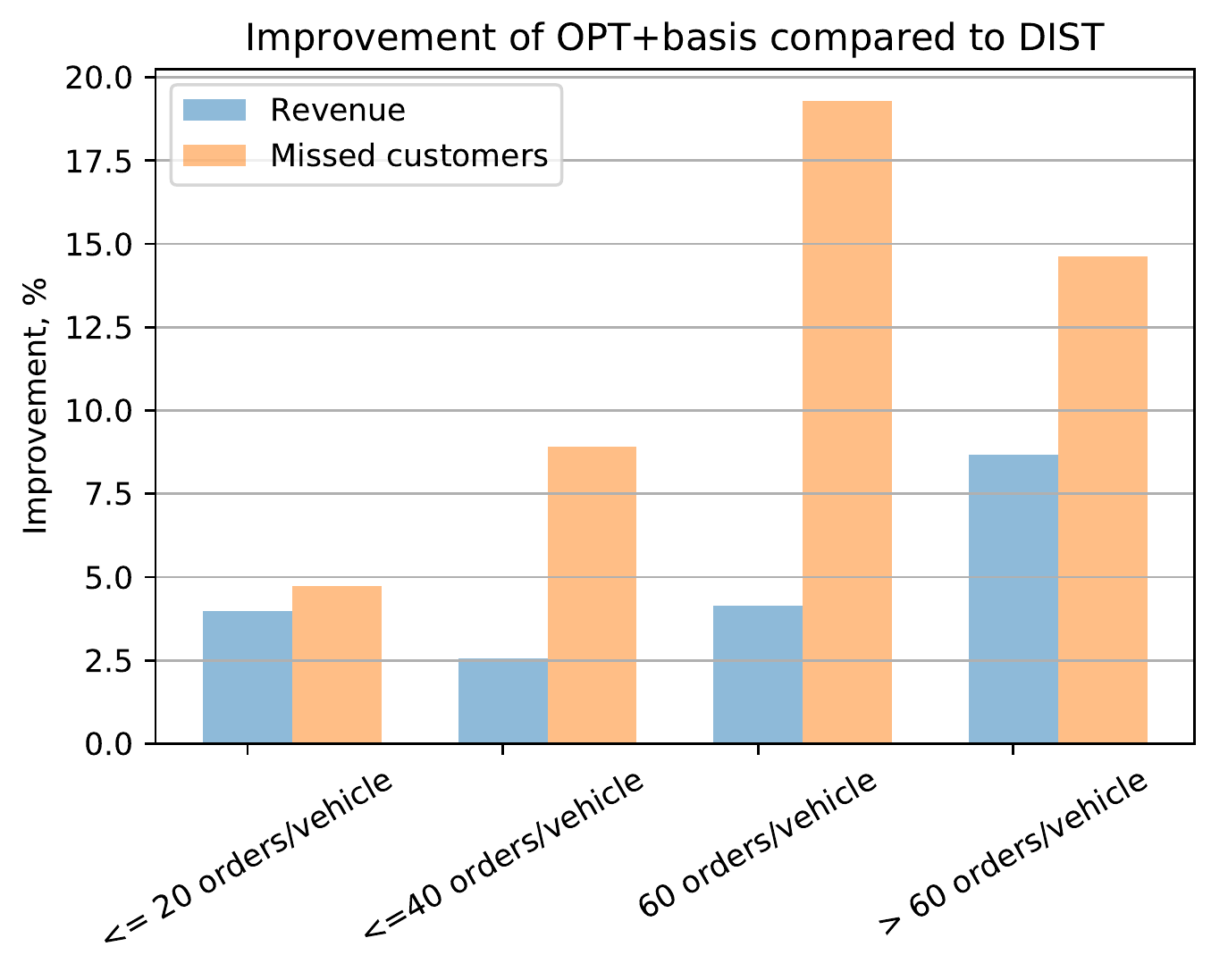}}
     \subfigure{\raisebox{13pt}{\includegraphics[width=0.47\columnwidth]{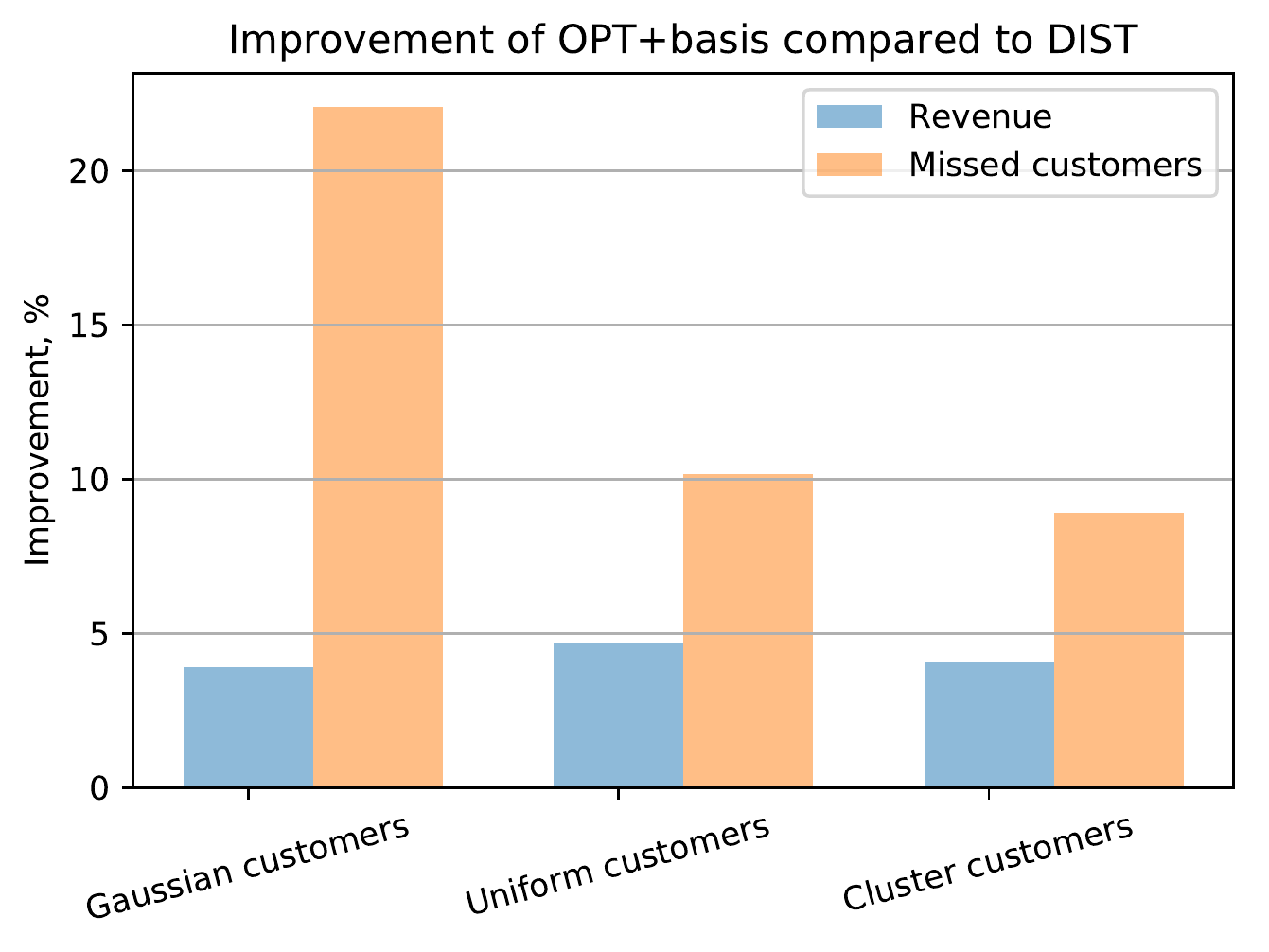}}}
    \caption{Percentage improvement in expected revenue and the reduction in missed customers between DIST and OPT+basis policies across different instance parameters.}
    \label{res_pri}
    \end{figure}
    
    \begin{figure}[h!]
    \centering
     \subfigure{\includegraphics[width=0.47\columnwidth]{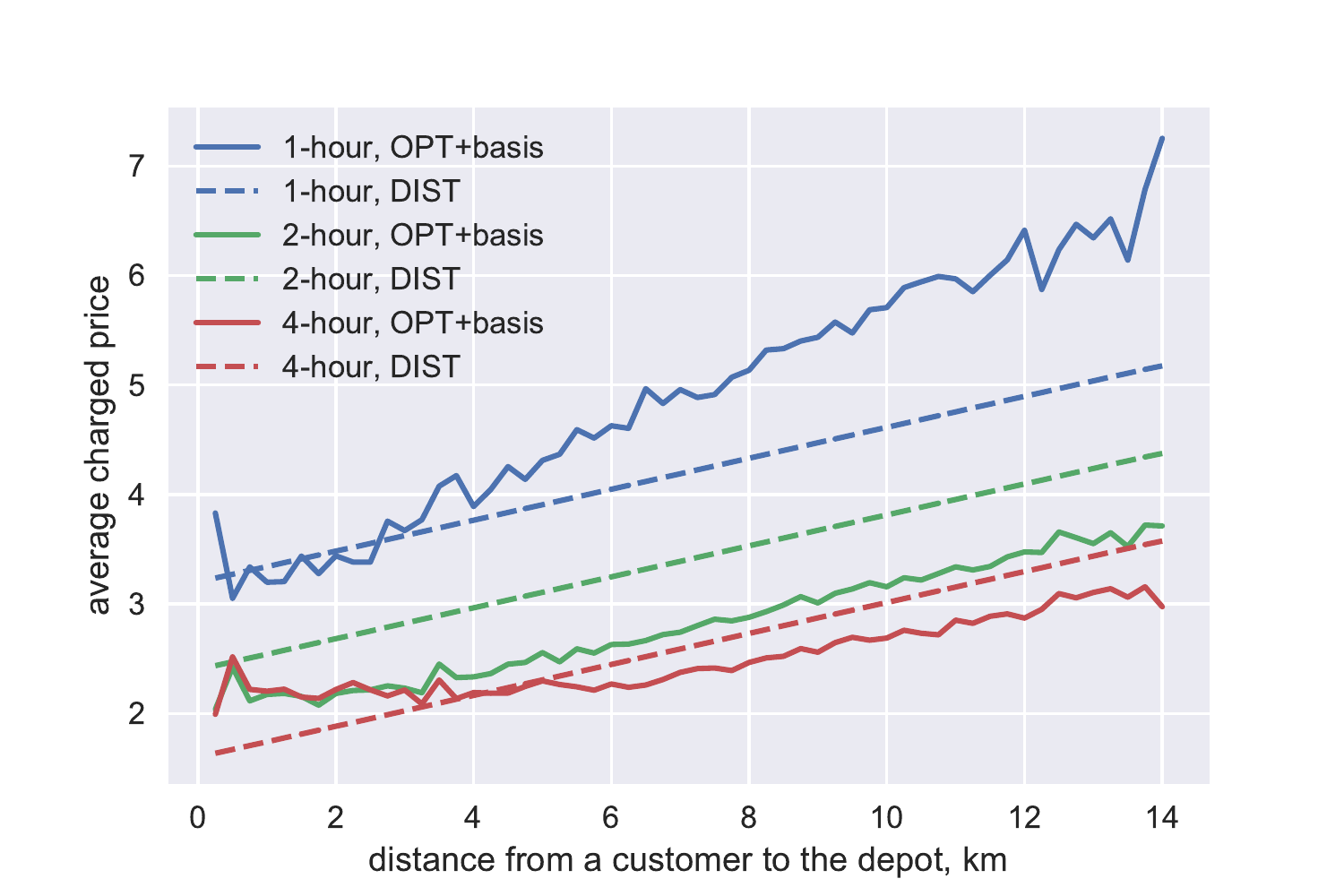}}
     \subfigure{\includegraphics[width=0.47\columnwidth]{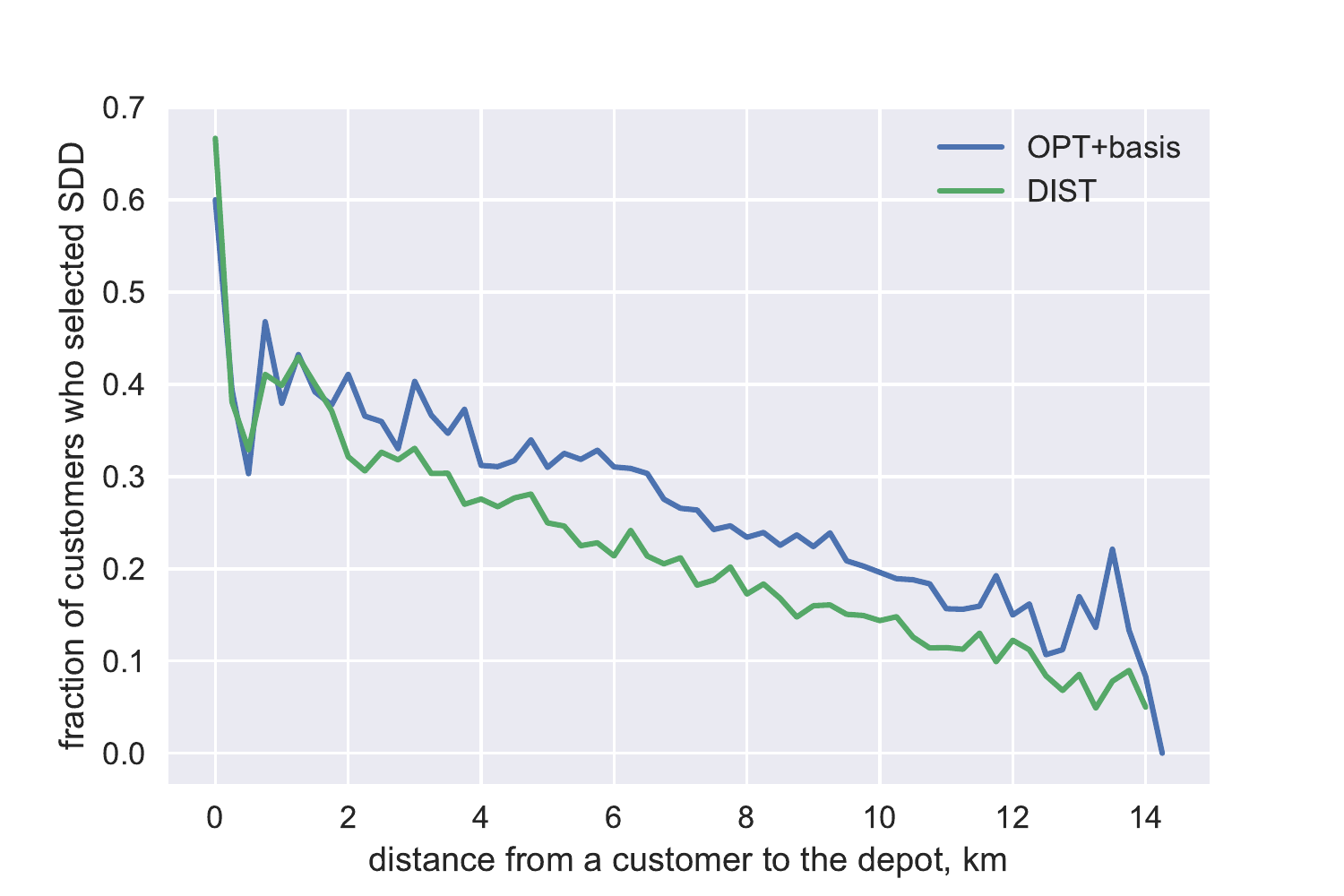}}
     \subfigure{\includegraphics[width=0.47\columnwidth]{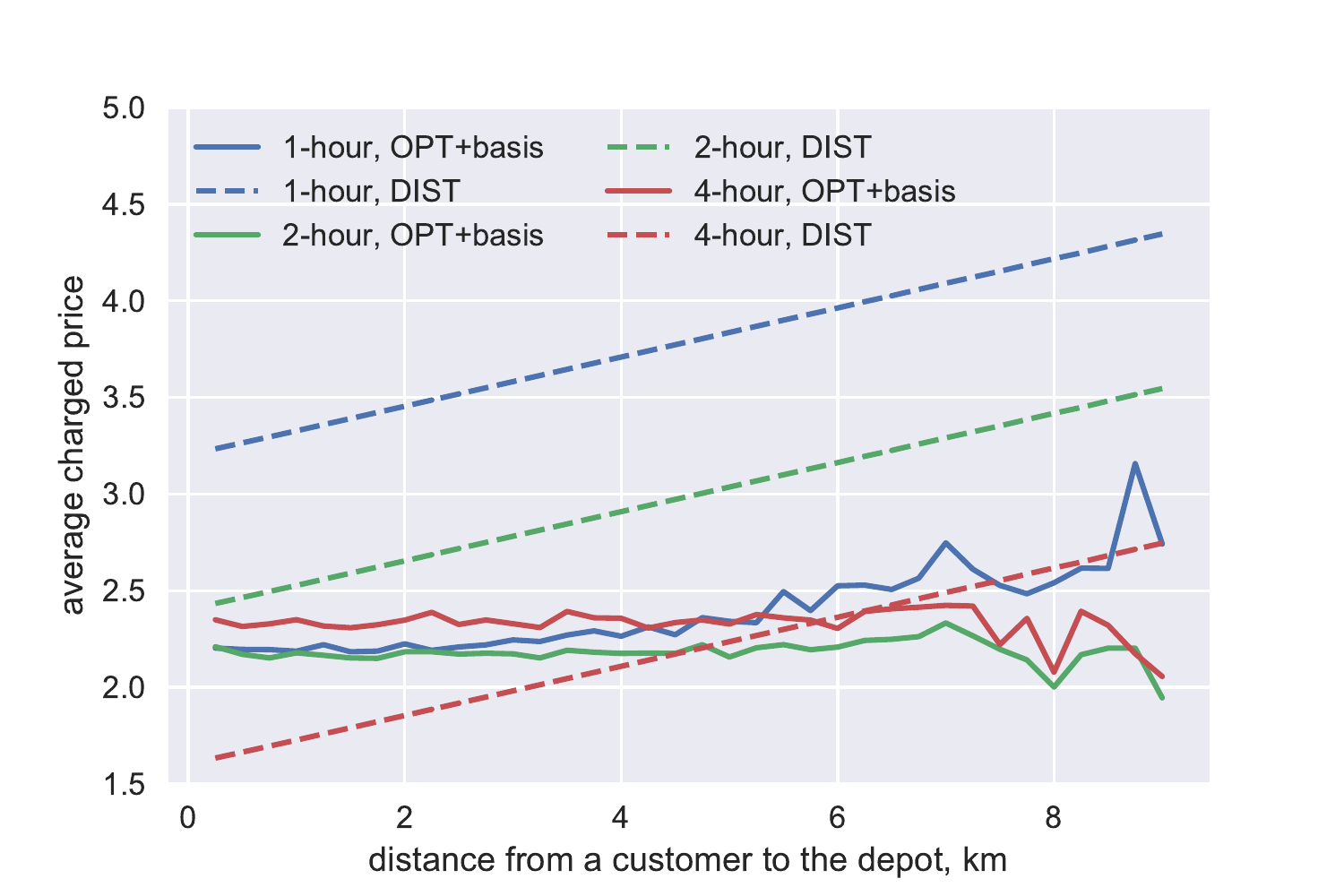}}
     \subfigure{\includegraphics[width=0.47\columnwidth]{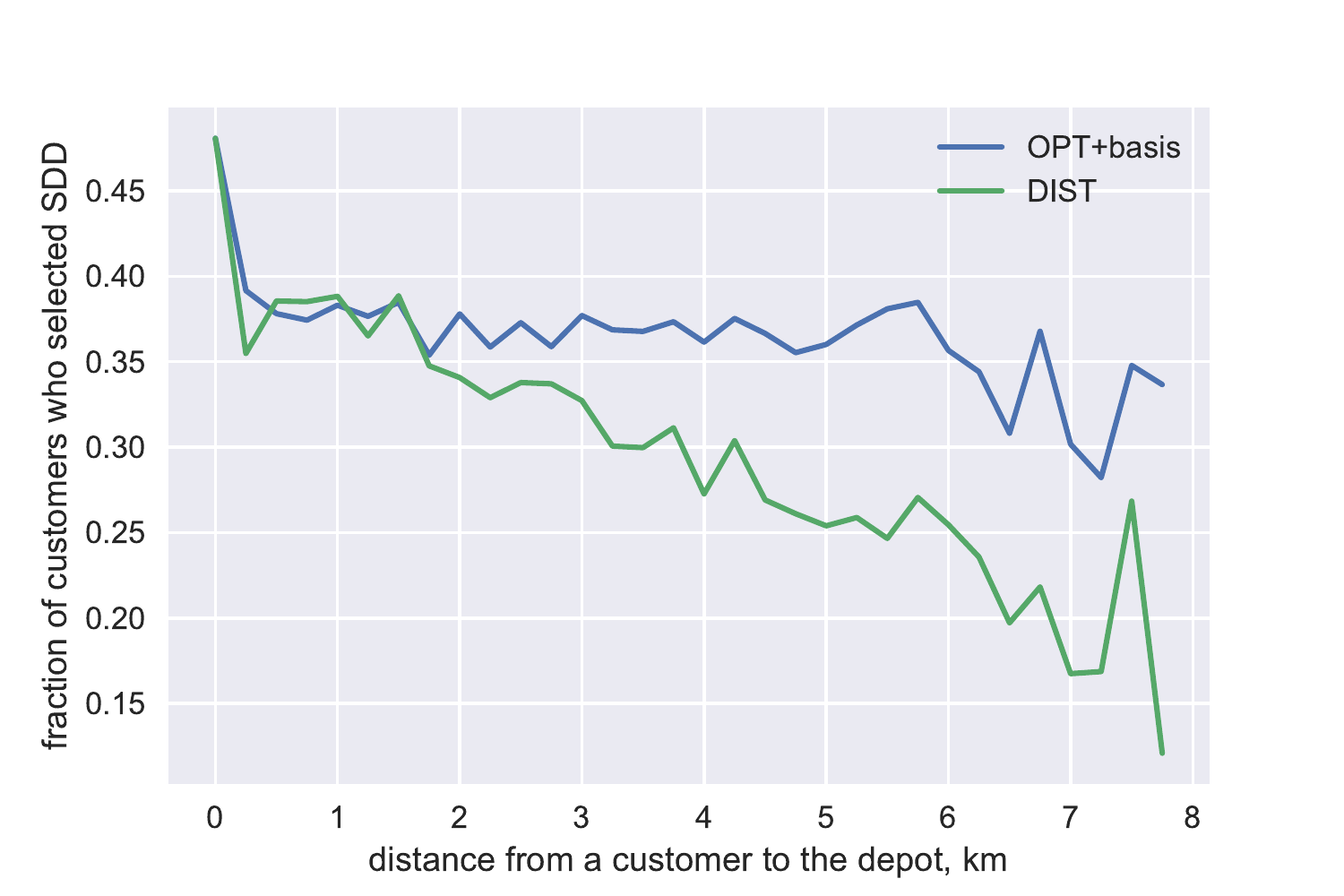}}

    \caption{Charged prices and proportion of customers selecting SDD for DIST and OPT+basis policies. Top: 80 uniform customers and 1 vehicle; bottom: 80 Gaussian customers and 1 vehicle. }
    \label{res_fairness}
    \end{figure}

\subsection{Computational results}
Here we conduct various experiments to explore how different routing and pricing methods influence the model performance.
\subsubsection{Performance Metrics} 
First we investigate appropriate metrics, beyond profit, to evaluate the different models. As any other companies, delivery companies aim to keep the customer satisfaction at a high level. Accordingly, we assess several customer-centric metrics: the number of customers whose requests for SSD were accepted (`Accepted'), the number of customers who were served within the deadline (`Served'), the number of customers whose delivery missed the deadline (`Missed'), the number of customers who rejected the provided SDD options and chose the next-day delivery (`Declined') and the number of customers for whom we were not able to provide SDD options (`Rejected'). Results are reported as averages of 1,000 runs.

\subsubsection{Stochastic routing}
One of the goals of this study is to exploit how modeling the stochasticity of travel times may lead to better pricing and routing decisions. Therefore, we are interested in how incorporating travel time distribution affects the model performance.
We compare three routing models: in one scenario only the expected travel times are taken into account, and travel times are hence treated in a deterministic manner (`Deterministic'),  in the second scenario the travel time distribution is fully known (`Stochastic'), and in the third scenario, the distributions are incorrectly specified. In the latter case, the true distribution might be Gaussian, while it is modeled as a mixture and vice versa (`Misspecified'). 
\par 
The probability of satisfying the deadline in the deterministic case is defined as:
\begin{equation}
     \mathbb{P}_R^{\text{const}}(d_{i_1}, \ldots d_{i_l}, t_{\text{deadline}}) = I\left\{ \sum_{j=1}^l d_{i_j}  < t_{\text{deadline}} \mu \right\},
\end{equation}
where $\{ d_{i_1}, \ldots d_{i_l} \}$ is the set of distances the vehicle has to traverse before the deadline $t_{\text{deadline}}$, and $\mu$ is the expected speed value. 

\subsubsection{Pricing policies} 
We compare the proposed VFA-based policy with fixed prices policy (where prices for each delivery deadline are predetermined) and distance-based pricing (where prices are proportionate to the distance between the customer and the depot). These baseline policies were also investigated in \cite{ulmer2017dynamic}. The baseline pricing is represented as: 
\begin{equation}
    \centering
    \mathcal{P}^{\text{FIX}} = [\alpha, 0.75\alpha, 0.5\alpha ],
    \end{equation}
\begin{equation}
\centering
    \mathcal{P}^{\text{DIST}} = [ \alpha + \gamma  d^i/d^{\text{max}}, 0.75  \alpha + \gamma  d^i/d^{\text{max}}, 0.5  \alpha + \gamma  d^i/d^{\text{max}}],
\end{equation}
where $\alpha$ and $\gamma$ are tunable parameters, $d^i$ is the distance between the current $i$-th customer and the depot, and $d^{\text{max}}$ is the maximum possible distance between any customer and the depot. We also compare the proposed approach with the policy proposed by \cite{ulmer2017dynamic}:
\begin{equation}
    P^i_k(C) = \max(p^i_b, \mathcal{O}^i(S_k)),
\end{equation} where $p^i_b$ are basis prices and  $\mathcal{O}^i(S_k)$ are the opportunity costs for the option $i$. We refer to this pricing policy as `OPP'. In all cases, the tunable parameters are found by the policy search approach. The proposed method is referred to as `OPT'. Table~\ref{tablefinalresall} shows the results for all combinations of pricing policies and travel time assumptions.

\begin{table}[htb]
\renewcommand{\arraystretch}{1.5}
\small
\begin{tabular}{|p{20mm}|p{20mm}|p{24mm}|p{24mm}|p{24mm}|p{24mm}|p{24mm}|}
\hline
\begin{tabular}[c]{@{}l@{}}Customer \\  Distribution\end{tabular} & \begin{tabular}[c]{@{}l@{}}Travel Time \\  assumption\end{tabular} & FIX       & DIST      & OPP       & OPT       & OPT+basis \\ \hline \hline
Gaussian                                                            & Deterministic                                                        & 61.9 / 32.3 / 1.4 & 63.1 / 32.2 / 1.7 & 61.3 / 32.3 / 1.7 & 61.4 / 35.6 / 3.8 & 64.6 / 31.9 / 1.9 \\ \hline
Gaussian                                                            & Stochastic                                                           & 61.2 / 32.9 / 1.5 & 63.1 / 34.0 / 1.7 & 61.2 / 34.4 / 1.3 & 64.8 / 34.8 / 1.6 & 65.5 / 33.9 / 1.3 \\ \hline
Gaussian                                                            & Misspecified                                             & 63.9 / 34.6 / 1.6 & 65.6 / 34.9 / 1.5 & 63.6 / 34.4 / 1.4 & 65.3 / 36.2 / 2.6 & 67.2 / 33.0 / 1.4 \\ \hline \hline
Uniform                                                             & Deterministic                                                        & 43.7 / 21.3 / 2.9 & 46.4 / 22.9 / 3.4 & 43.0 / 22.7 / 3.5 & 33.3 / 31.3 / 10 & 46.8 / 22.5 / 3.1 \\ \hline
Uniform                                                             & Stochastic                                                           & 44.2 / 23.0 / 3.2 & 46.8 / 22.6 / 2.8 & 45.7 / 23.6 / 2.9 & 41.4 / 30.9 / 6.8 & 49.0 / 23.4 / 2.5 \\ \hline
Uniform                                                             & Misspecified                                            & 45.0 / 22.7 / 3.3 & 46.5 / 24.4 / 3.8 & 44.9 / 23.7 / 3.4 & 37.9 / 31.2 / 8.2 & 47.8 / 22.8 / 2.7 \\ \hline \hline
Cluster                                                             & Deterministic                                                        & 46.6 / 21.9 / 2.4 & 48.8 / 22.5 / 2.9 & 45.8 / 23.4 / 3.0 & 37.8 / 32.6 / 9.4 & 48.1 / 23.0 / 3.1 \\ \hline
Cluster                                                             & Stochastic                                                           & 48.0 / 23.4 / 2.4 & 49.6 / 22.8 / 2.7 & 49.4 / 24.4 / 2.1 & 47.2 / 32.2 / 5.5 & 51.6 / 24.8 / 2.5 \\ \hline
Cluster                                                             & Misspecified                                             & 48.7 / 24.4 / 2.6 & 49.3 / 22.9 / 2.6 & 48.4 / 23.9 / 2.5 & 43.8 / 32.5 / 7.1 & 50.5 / 24.9 / 2.8 \\ \hline
\end{tabular}
\caption{Comparison between models with different TT assumptions and pricing policies for Gaussian travel times. Values in cells are average revenue/ average number of accepted customers/ average number of missed customers. Values are averaged over all test instance specifications.}
\label{tablefinalresall}
\end{table}

\subsubsection{Basis Prices}
As mentioned in the previous section, \cite{ulmer2017dynamic} proposes a pricing policy with basis prices. In that paper, it is shown that without basis prices, VFA is unable to correctly approximate values. This approach is similar to certain services like Uber, where in cases of high supply basis prices are offered, and in cases with high opportunity costs and/or high demand, the prices are increased (surge pricing). \par
Similarly, we perform experiments with basis prices and compare the results. In the proposed approach, basis prices are set as a lower bound for L-BFGS optimization. To make this policy more competitive, the basis prices are proportional to the distance from the depot (similar to DIST policy). This policy is referred to as `OPT+basis' in Table~\ref{tablefinalresall}. \par
It should be noted that in all policies there are no restrictions on the ordering of prices (i.e. the price of 4-hour delivery can be higher than the price of 1-hour delivery).

\subsection{Analysis of results}

The analysis in this section has two objectives: 1) to investigate whether the travel time distribution information is crucial for same-day delivery routing and pricing; 2) to benchmark the proposed approach for pricing of SDD. Table~\ref{tablefinalresall} provides an overview of results that are averaged across several instances with varying numbers of orders (40, 80 and 120), vehicles (1 to 3), and different penalties (0 and 2). We will discuss those results in more detail in this section. 
\par 
Figure~\ref{res_imp} shows that the proposed model outperforms the baseline with deterministic travel times in terms of revenue and served customers. However, the improvement is only marginal when there are only few orders (not more than 20 orders per vehicle) or customer locations are distributed according to Gaussian distribution. Conversely, the model which incorporates the travel time distributions performs best in cases when the demand is higher than the supply or when the relative distances between customers and depot are large (customers locations are uniform or cluster distributed). 
\par 
Figure~\ref{res_pri} compares the performance of the proposed OPT+basis policy and the baseline DIST policy. The largest profit gain is observed for instances with larger numbers of orders per vehicle. For such instances the reduction in the number of missed customers is also more pronounced.
\par 
Table~\ref{simple_comp} shows the distribution of selected delivery options and the average price paid for eight different instances. We are interested in how the supply and demand levels affect both delivery prices and SDD service levels. For this reason, the table contains instances with low supply (80 orders and 1 vehicle), high supply (80 orders and 3 vehicles), as well as with customer locations distributed closer to the depot (Gaussian) and further away (Uniform). It can be observed that in cases with low supply (1 vehicle per 80 orders) the policy rejects up to 20\% of the customers by not providing any same-day delivery options. However, in cases with more supply (3 vehicles per 80 orders) virtually all customers have an opportunity to choose SDD. In cases when a penalty for missed deliveries is present, the delivery prices are almost always higher since they have to absorb possible penalties as well as control the number of customers that select same-day delivery. The difference in prices between the penalty and non-penalty situations is significantly higher in low supply cases.
\par
Full results are presented in the appendix.

\subsection{Fairness}
\cite{fairness} discuss the issue of fairness with regards to the customer accepting and pricing mechanisms in vehicle routing problems. They show that a policy that can decline certain customer requests can overall accept more requests. At the same time, this policy underserves areas distant from the depot. We aim to investigate this issue in more detail. To this end, we compare the distance-based baseline policy (DIST) with our proposed method (OPT+basis). DIST policy by design discriminates against customers farther from the depot by increasing prices proportionally to the distance. Figure~\ref{res_fairness} provides a comparison between the two policies in terms of both prices and SDD services. For the instance with uniform customers, both policies lead to a smaller proportion of distant customers selecting SDD. However, the opportunity costs based approach achieves higher SDD service rates compared to DIST policy. As can be seen on the price graph, OPT+basis prices 1-hour option significantly higher than DIST, while 2- and 4-hour options are cheaper with only a small increase in price according to distance. This means that for OPT+basis `distant' customers realistically only have 2- and 4-hour options but they are competitively priced. For DIST, the distant customers can select any option but all of them are priced higher than average. For the second instance, where the customers are Gaussian distributed in space, the difference between the two policies is even more pronounced. OPT+basis policy provides an almost constant rate of SDD deliveries across the whole service area. The prices charged are similar for most of the area as well. One of the possible reasons is that in this case, the number of distant customers is relatively low (e.g. any coordinate that is more than 5 km away from the depot already represents 2 standard deviations of the spatial distribution). This illustrates a tradeoff between the number of accepted requests and the average price of the delivery. In real-life situations, companies should either choose pricing policies which do not overprice certain parts of the service area or reduce the area altogether.

\par 

\begin{figure}[htb]
    \centering
     \subfigure{\includegraphics[width=0.67\columnwidth]{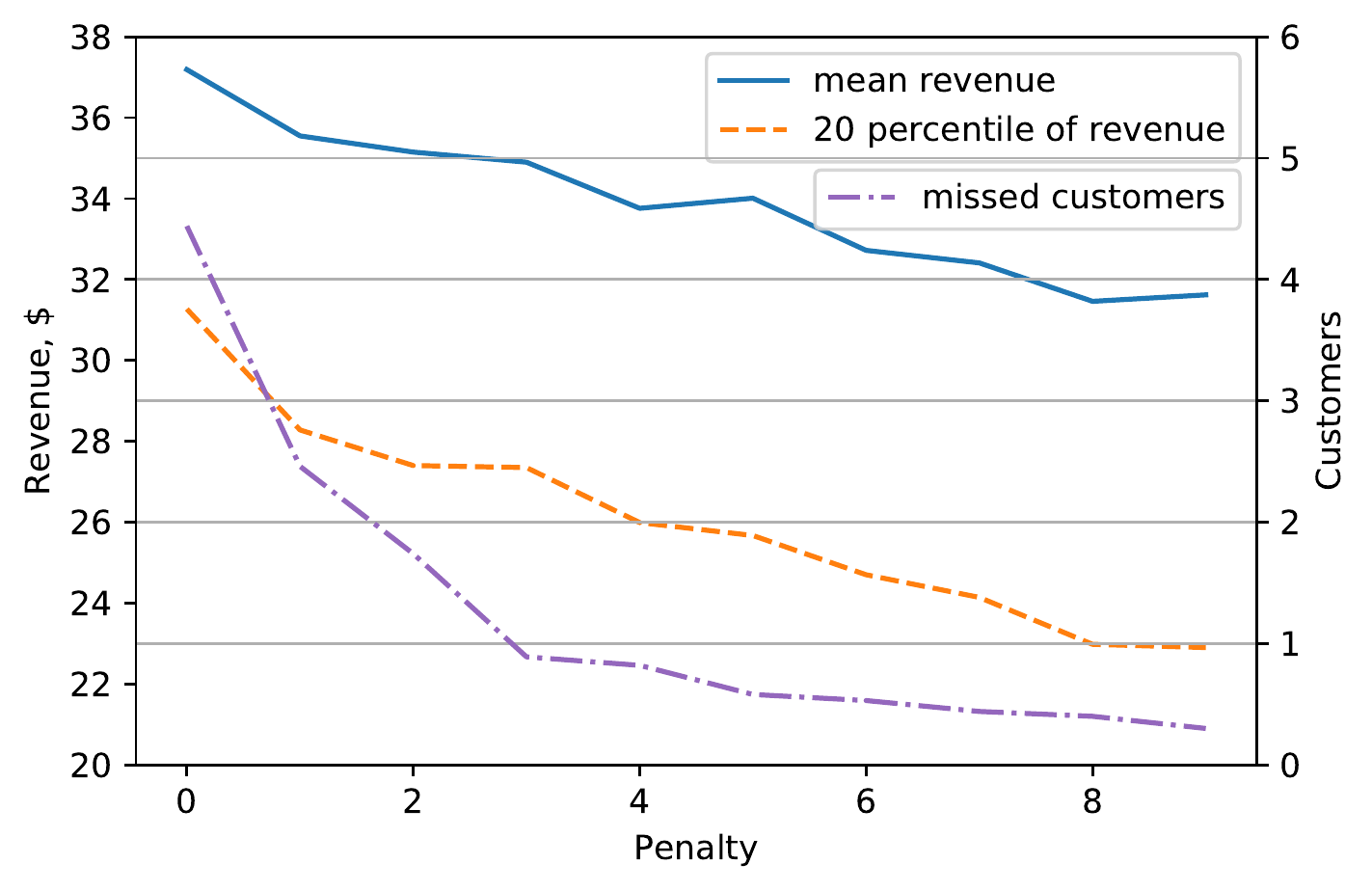}}
     \caption{The effect of the penalty for missed deliveries (Instance parameters: 80 orders, 1 vehicle, Gaussian model, Gaussian travel times, Gaussian customer distribution).}
    \label{figpen}
    \end{figure}

\section{Conclusion}
In this study, we have proposed a method for simultaneously optimizing routing and pricing for same-day delivery routing that also takes into account the variability of travel times. We have utilized value function approximation to compute opportunity costs for accepting a given customer request. We have performed an extensive set of simulations and compared the proposed method with conventional pricing policies and a model with deterministic travel times. We have shown that information about the travel time distribution can greatly improve the quality of routing and pricing solution. We have also investigated several new issues arising due to the stochastic nature of same-day delivery. Specifically, we have investigated how penalties for late deliveries affect the pricing structure. Moreover, we have analyzed the trade-off between the revenue and missed deliveries as well as the issue of the pricing fairness for different policies. 
\par 
We will further explore several aspects of this problem in future work. First, the model can be extended to support incomplete information about travel time distributions (e.g., the distribution belongs to a family of distributions). Next, the model can be improved by employing different routing heuristics, as well as other approaches for approximating opportunity costs. Additionally, the model can be modified to incorporate ad-hoc delivery drivers, in which case the drivers can choose orders depending on the compensation.
\par 
In this study, we have considered linear regression as an estimator for VFA. However, for real-world applications, the prediction performance of linear models can be insufficient, especially since the set of features (`basis functions') have to be manually constructed. In some cases, a model-free reinforcement learning approach can be more beneficial. A book by \cite{busoniu2017reinforcement} provides a thorough overview of both (approximate) dynamic programming and reinforcement learning. 

\section*{Acknowledgement}
This work was partially supported by the Singapore National Research Foundation through the Singapore-MIT Alliance for Research and Technology (SMART) Centre for Future Urban Mobility (FM).

\newcommand{\newblock}{}
\bibliographystyle{apacite}
\bibliography{sample}

\begin{APPENDICES}
\section{Full results}
The appendix contains full results for each instance specification (the number of orders, the number of vehicles, and the penalty for missed delivery) and each model specification (the customer distribution, the travel time distribution assumption, and the pricing policy). All instances are described in Table~\ref{appendix_instance}. The results in Tables~\ref{first_full_table}-\ref{last_full_table} are averaged over 1000 simulation runs.
\begin{table}[htb]
            \renewcommand{\arraystretch}{1.0}
            \centering
\begin{tabular}{|l|l|l|l|}
            \hline
            Instance Id & Orders & Vehicles & Penalty \\ \hline
0 & 40 & 1 & 0 \\ \hline
1 & 40 & 1 & 2 \\ \hline
2 & 40 & 2 & 0 \\ \hline
3 & 40 & 2 & 2 \\ \hline
4 & 40 & 3 & 0 \\ \hline
5 & 40 & 3 & 2 \\ \hline
6 & 80 & 1 & 0 \\ \hline
7 & 80 & 1 & 2 \\ \hline
8 & 80 & 2 & 0 \\ \hline
9 & 80 & 2 & 2 \\ \hline
10 & 80 & 3 & 0 \\ \hline
11 & 80 & 3 & 2 \\ \hline
12 & 120 & 1 & 0 \\ \hline
13 & 120 & 1 & 2 \\ \hline
14 & 120 & 2 & 0 \\ \hline
15 & 120 & 2 & 2 \\ \hline
16 & 120 & 3 & 0 \\ \hline
17 & 120 & 3 & 2 \\ \hline

\end{tabular}
\caption{Instance specification.} \label{appendix_instance} 
\end{table}

\begin{table}
\parbox{.45\linewidth}{
\centering
\begin{tabular}{p{5mm}p{14mm}p{10mm}p{15mm}p{11mm}}
            \hline
            Id & Revenue & Served & Accepted & Missed \\ \hline
            
0 & 35.32 & 21.44 & 21.82 & 0.38 \\ \hline 
1 & 35.17 & 18.84 & 19.02 & 0.18 \\ \hline 
2 & 36.14 & 21.54 & 21.78 & 0.24 \\ \hline 
3 & 35.65 & 21.58 & 21.83 & 0.25 \\ \hline 
4 & 36.1 & 21.48 & 21.74 & 0.26 \\ \hline 
5 & 35.57 & 21.5 & 21.78 & 0.28 \\ \hline 
6 & 60.0 & 26.01 & 27.47 & 1.46 \\ \hline 
7 & 57.45 & 26.04 & 27.39 & 1.35 \\ \hline 
8 & 69.42 & 33.97 & 35.03 & 1.06 \\ \hline 
9 & 67.49 & 36.39 & 37.71 & 1.32 \\ \hline 
10 & 71.29 & 36.55 & 37.68 & 1.13 \\ \hline 
11 & 68.98 & 36.53 & 37.71 & 1.18 \\ \hline 
12 & 71.07 & 27.4 & 29.93 & 2.53 \\ \hline 
13 & 67.03 & 25.88 & 27.3 & 1.42 \\ \hline 
14 & 90.42 & 40.2 & 42.85 & 2.65 \\ \hline 
15 & 85.33 & 40.16 & 42.72 & 2.56 \\ \hline 
16 & 99.45 & 50.66 & 54.13 & 3.47 \\ \hline 
17 & 92.42 & 50.56 & 54.05 & 3.49
\\ \hline 

\end{tabular}
\caption{Result of FIX strategy on Gaussian distributed customers with Deterministic travel time distribution assumption.}
\label{first_full_table}
}
\hfill
\parbox{.45\linewidth}{
\centering
\begin{tabular}{p{5mm}p{14mm}p{10mm}p{15mm}p{11mm}}
            \hline
            Id & Revenue & Served & Accepted & Missed \\ \hline
            
0 & 27.83 & 12.65 & 14.85 & 2.2 \\ \hline 
1 & 25.33 & 10.73 & 11.48 & 0.75 \\ \hline 
2 & 33.48 & 15.78 & 16.18 & 0.4 \\ \hline 
3 & 32.83 & 16.81 & 17.34 & 0.53 \\ \hline 
4 & 35.21 & 18.44 & 18.72 & 0.28 \\ \hline 
5 & 34.63 & 18.4 & 18.67 & 0.27 \\ \hline 
6 & 38.06 & 12.69 & 14.26 & 1.57 \\ \hline 
7 & 33.98 & 13.35 & 15.72 & 2.37 \\ \hline 
8 & 54.59 & 23.33 & 27.4 & 4.07 \\ \hline 
9 & 48.67 & 19.53 & 20.91 & 1.38 \\ \hline 
10 & 62.46 & 27.64 & 29.73 & 2.09 \\ \hline 
11 & 58.57 & 25.94 & 27.39 & 1.45 \\ \hline 
12 & 44.33 & 14.3 & 18.34 & 4.04 \\ \hline 
13 & 38.98 & 12.3 & 13.46 & 1.16 \\ \hline 
14 & 66.31 & 26.83 & 35.91 & 9.08 \\ \hline 
15 & 57.52 & 22.19 & 24.81 & 2.62 \\ \hline 
16 & 76.75 & 31.02 & 35.93 & 4.91 \\ \hline 
17 & 68.6 & 29.39 & 32.82 & 3.43
\\ \hline 

\end{tabular}
\caption{Result of FIX strategy on Cluster distributed customers with Deterministic travel time distribution assumption.}
}
\end{table}

\begin{table}
\parbox{.45\linewidth}{
\centering
\begin{tabular}{p{5mm}p{14mm}p{10mm}p{15mm}p{11mm}}
            \hline
            Id & Revenue & Served & Accepted & Missed \\ \hline
            
0 & 26.49 & 11.02 & 12.71 & 1.69 \\ \hline 
1 & 23.66 & 10.42 & 11.62 & 1.2 \\ \hline 
2 & 32.68 & 15.45 & 16.29 & 0.84 \\ \hline 
3 & 31.29 & 14.52 & 15.04 & 0.52 \\ \hline 
4 & 35.25 & 18.41 & 18.87 & 0.46 \\ \hline 
5 & 34.29 & 18.39 & 18.92 & 0.53 \\ \hline 
6 & 35.79 & 12.31 & 15.68 & 3.37 \\ \hline 
7 & 28.98 & 12.3 & 15.64 & 3.34 \\ \hline 
8 & 51.29 & 20.96 & 24.96 & 4.0 \\ \hline 
9 & 44.53 & 21.26 & 24.93 & 3.67 \\ \hline 
10 & 59.94 & 28.0 & 32.11 & 4.11 \\ \hline 
11 & 55.18 & 23.58 & 24.99 & 1.41 \\ \hline 
12 & 39.89 & 12.49 & 16.79 & 4.3 \\ \hline 
13 & 35.41 & 11.73 & 13.66 & 1.93 \\ \hline 
14 & 61.05 & 22.92 & 30.15 & 7.23 \\ \hline 
15 & 52.25 & 21.13 & 24.98 & 3.85 \\ \hline 
16 & 73.48 & 29.62 & 35.97 & 6.35 \\ \hline 
17 & 64.9 & 26.81 & 30.12 & 3.31
\\ \hline

\end{tabular}
\caption{Result of FIX strategy on Uniform distributed customers with Deterministic travel time distribution assumption.}
}
\hfill
\parbox{.45\linewidth}{
\centering
\begin{tabular}{p{5mm}p{14mm}p{10mm}p{15mm}p{11mm}}
            \hline
            Id & Revenue & Served & Accepted & Missed \\ \hline
            
0 & 34.7 & 18.48 & 19.0 & 0.52 \\ \hline 
1 & 34.08 & 17.39 & 17.7 & 0.31 \\ \hline 
2 & 36.14 & 21.57 & 21.78 & 0.21 \\ \hline 
3 & 35.71 & 21.58 & 21.78 & 0.2 \\ \hline 
4 & 36.18 & 21.56 & 21.82 & 0.26 \\ \hline 
5 & 35.61 & 21.53 & 21.79 & 0.26 \\ \hline 
6 & 55.33 & 23.06 & 25.14 & 2.08 \\ \hline 
7 & 51.45 & 23.08 & 25.16 & 2.08 \\ \hline 
8 & 70.77 & 36.8 & 37.66 & 0.86 \\ \hline 
9 & 69.24 & 36.87 & 37.71 & 0.84 \\ \hline 
10 & 71.33 & 36.54 & 37.67 & 1.13 \\ \hline 
11 & 69.18 & 36.6 & 37.75 & 1.15 \\ \hline 
12 & 64.28 & 24.09 & 27.25 & 3.16 \\ \hline 
13 & 57.26 & 23.78 & 27.25 & 3.47 \\ \hline 
14 & 94.12 & 44.06 & 46.45 & 2.39 \\ \hline 
15 & 91.32 & 41.54 & 42.81 & 1.27 \\ \hline 
16 & 100.09 & 57.86 & 61.88 & 4.02 \\ \hline 
17 & 94.05 & 58.29 & 61.89 & 3.6
\\ \hline 

\end{tabular}
\caption{Result of FIX strategy on Gaussian distributed customers with Gaussian travel time distribution assumption.}
}
\end{table}

\begin{table}
\parbox{.45\linewidth}{
\centering
\begin{tabular}{p{5mm}p{14mm}p{10mm}p{15mm}p{11mm}}
            \hline
            Id & Revenue & Served & Accepted & Missed \\ \hline
            
0 & 26.67 & 11.16 & 12.51 & 1.35 \\ \hline 
1 & 23.89 & 11.11 & 12.52 & 1.41 \\ \hline 
2 & 34.33 & 18.12 & 18.7 & 0.58 \\ \hline 
3 & 32.46 & 14.7 & 14.9 & 0.2 \\ \hline 
4 & 35.33 & 18.44 & 18.71 & 0.27 \\ \hline 
5 & 34.91 & 18.47 & 18.72 & 0.25 \\ \hline 
6 & 37.66 & 13.16 & 15.77 & 2.61 \\ \hline 
7 & 33.16 & 12.33 & 14.29 & 1.96 \\ \hline 
8 & 55.33 & 22.81 & 25.07 & 2.26 \\ \hline 
9 & 53.91 & 23.42 & 25.17 & 1.75 \\ \hline 
10 & 64.36 & 29.97 & 32.4 & 2.43 \\ \hline 
11 & 62.01 & 28.52 & 29.86 & 1.34 \\ \hline 
12 & 41.08 & 13.64 & 18.36 & 4.72 \\ \hline 
13 & 37.68 & 12.03 & 13.45 & 1.42 \\ \hline 
14 & 65.67 & 29.27 & 39.11 & 9.84 \\ \hline 
15 & 66.92 & 29.2 & 32.93 & 3.73 \\ \hline 
16 & 78.62 & 32.08 & 35.95 & 3.87 \\ \hline 
17 & 79.44 & 38.55 & 42.64 & 4.09
\\ \hline

\end{tabular}
\caption{Result of FIX strategy on Cluster distributed customers with Gaussian travel time distribution assumption.}
}
\hfill
\parbox{.45\linewidth}{
\centering
\begin{tabular}{p{5mm}p{14mm}p{10mm}p{15mm}p{11mm}}
            \hline
            Id & Revenue & Served & Accepted & Missed \\ \hline
            
0 & 24.77 & 10.49 & 12.69 & 2.2 \\ \hline 
1 & 22.67 & 10.17 & 11.61 & 1.44 \\ \hline 
2 & 33.14 & 15.68 & 16.33 & 0.65 \\ \hline 
3 & 32.27 & 16.82 & 17.56 & 0.74 \\ \hline 
4 & 35.34 & 18.45 & 18.9 & 0.45 \\ \hline 
5 & 34.41 & 18.46 & 18.95 & 0.49 \\ \hline 
6 & 33.2 & 10.9 & 12.86 & 1.96 \\ \hline 
7 & 29.77 & 9.59 & 10.4 & 0.81 \\ \hline 
8 & 50.61 & 23.23 & 29.63 & 6.4 \\ \hline 
9 & 48.14 & 23.6 & 27.22 & 3.62 \\ \hline 
10 & 58.84 & 31.04 & 37.4 & 6.36 \\ \hline 
11 & 56.44 & 31.12 & 34.73 & 3.61 \\ \hline 
12 & 38.54 & 12.05 & 15.07 & 3.02 \\ \hline 
13 & 33.72 & 10.05 & 11.0 & 0.95 \\ \hline 
14 & 59.35 & 25.26 & 35.89 & 10.63 \\ \hline 
15 & 60.45 & 22.76 & 24.95 & 2.19 \\ \hline 
16 & 72.08 & 31.2 & 39.3 & 8.1 \\ \hline 
17 & 72.44 & 34.52 & 39.26 & 4.74
\\ \hline

\end{tabular}
\caption{Result of FIX strategy on Uniform distributed customers with Gaussian travel time distribution assumption.}
}
\end{table}

\begin{table}
\parbox{.45\linewidth}{
\centering
\begin{tabular}{p{5mm}p{14mm}p{10mm}p{15mm}p{11mm}}
            \hline
            Id & Revenue & Served & Accepted & Missed \\ \hline
            
0 & 35.79 & 21.61 & 21.73 & 0.12 \\ \hline 
1 & 35.51 & 21.57 & 21.71 & 0.14 \\ \hline 
2 & 35.89 & 18.79 & 18.95 & 0.16 \\ \hline 
3 & 35.56 & 21.45 & 21.7 & 0.25 \\ \hline 
4 & 36.1 & 21.49 & 21.75 & 0.26 \\ \hline 
5 & 35.49 & 21.42 & 21.68 & 0.26 \\ \hline 
6 & 63.74 & 30.43 & 32.42 & 1.99 \\ \hline 
7 & 60.12 & 28.44 & 29.9 & 1.46 \\ \hline 
8 & 69.95 & 41.31 & 43.06 & 1.75 \\ \hline 
9 & 68.75 & 36.57 & 37.52 & 0.95 \\ \hline 
10 & 71.12 & 36.52 & 37.6 & 1.08 \\ \hline 
11 & 68.87 & 36.47 & 37.63 & 1.16 \\ \hline 
12 & 77.92 & 34.2 & 39.28 & 5.08 \\ \hline 
13 & 73.39 & 30.77 & 32.87 & 2.1 \\ \hline 
14 & 93.51 & 46.15 & 50.19 & 4.04 \\ \hline 
15 & 92.15 & 44.52 & 46.53 & 2.01 \\ \hline 
16 & 100.85 & 51.13 & 54.07 & 2.94 \\ \hline 
17 & 96.12 & 51.38 & 54.06 & 2.68
\\ \hline

\end{tabular}
\caption{Result of FIX strategy on Gaussian distributed customers with Misspecified travel time distribution assumption.}
}
\hfill
\parbox{.45\linewidth}{
\centering
\begin{tabular}{p{5mm}p{14mm}p{10mm}p{15mm}p{11mm}}
            \hline
            Id & Revenue & Served & Accepted & Missed \\ \hline
            
0 & 29.07 & 14.07 & 16.2 & 2.13 \\ \hline 
1 & 26.17 & 13.29 & 14.89 & 1.6 \\ \hline 
2 & 34.08 & 16.9 & 17.39 & 0.49 \\ \hline 
3 & 33.65 & 18.19 & 18.7 & 0.51 \\ \hline 
4 & 35.18 & 20.93 & 21.37 & 0.44 \\ \hline 
5 & 34.47 & 20.98 & 21.35 & 0.37 \\ \hline 
6 & 42.43 & 16.6 & 21.01 & 4.41 \\ \hline 
7 & 37.74 & 14.07 & 15.71 & 1.64 \\ \hline 
8 & 54.39 & 23.48 & 27.36 & 3.88 \\ \hline 
9 & 52.67 & 24.61 & 27.35 & 2.74 \\ \hline 
10 & 63.39 & 29.57 & 32.36 & 2.79 \\ \hline 
11 & 60.78 & 28.17 & 29.73 & 1.56 \\ \hline 
12 & 50.41 & 17.13 & 20.36 & 3.23 \\ \hline 
13 & 43.98 & 15.74 & 18.24 & 2.5 \\ \hline 
14 & 64.41 & 25.92 & 32.84 & 6.92 \\ \hline 
15 & 62.81 & 28.08 & 32.79 & 4.71 \\ \hline 
16 & 75.82 & 29.71 & 32.92 & 3.21 \\ \hline 
17 & 75.17 & 35.01 & 39.33 & 4.32
\\ \hline

\end{tabular}
\caption{Result of FIX strategy on Cluster distributed customers with Misspecified travel time distribution assumption.}
}
\end{table}

\begin{table}
\parbox{.45\linewidth}{
\centering
\begin{tabular}{p{5mm}p{14mm}p{10mm}p{15mm}p{11mm}}
            \hline
            Id & Revenue & Served & Accepted & Missed \\ \hline
            
0 & 27.0 & 12.47 & 15.05 & 2.58 \\ \hline 
1 & 23.26 & 11.84 & 13.86 & 2.02 \\ \hline 
2 & 33.2 & 17.54 & 18.87 & 1.33 \\ \hline 
3 & 31.76 & 17.81 & 18.89 & 1.08 \\ \hline 
4 & 35.39 & 18.45 & 18.88 & 0.43 \\ \hline 
5 & 34.14 & 20.97 & 21.53 & 0.56 \\ \hline 
6 & 38.05 & 14.35 & 18.96 & 4.61 \\ \hline 
7 & 32.29 & 12.96 & 15.68 & 2.72 \\ \hline 
8 & 49.91 & 21.62 & 27.18 & 5.56 \\ \hline 
9 & 48.63 & 20.72 & 22.86 & 2.14 \\ \hline 
10 & 60.05 & 26.72 & 29.59 & 2.87 \\ \hline 
11 & 57.31 & 25.56 & 27.23 & 1.67 \\ \hline 
12 & 43.85 & 15.59 & 22.67 & 7.08 \\ \hline 
13 & 39.67 & 13.26 & 15.16 & 1.9 \\ \hline 
14 & 56.16 & 24.38 & 39.35 & 14.97 \\ \hline 
15 & 58.63 & 22.34 & 25.05 & 2.71 \\ \hline 
16 & 71.06 & 26.84 & 30.22 & 3.38 \\ \hline 
17 & 69.04 & 25.93 & 27.41 & 1.48
\\ \hline

\end{tabular}
\caption{Result of FIX strategy on Uniform distributed customers with Misspecified travel time distribution assumption.}
}
\hfill
\parbox{.45\linewidth}{
\centering
\begin{tabular}{p{5mm}p{14mm}p{10mm}p{15mm}p{11mm}}
            \hline
            Id & Revenue & Served & Accepted & Missed \\ \hline
            
0 & 36.23 & 20.5 & 20.77 & 0.27 \\ \hline 
1 & 35.52 & 21.16 & 21.45 & 0.29 \\ \hline 
2 & 36.81 & 21.0 & 21.24 & 0.24 \\ \hline 
3 & 36.23 & 21.21 & 21.43 & 0.22 \\ \hline 
4 & 36.8 & 21.24 & 21.48 & 0.24 \\ \hline 
5 & 36.4 & 21.03 & 21.24 & 0.21 \\ \hline 
6 & 62.68 & 27.49 & 29.14 & 1.65 \\ \hline 
7 & 60.29 & 26.8 & 27.98 & 1.18 \\ \hline 
8 & 70.84 & 40.61 & 42.54 & 1.93 \\ \hline 
9 & 67.95 & 38.06 & 39.61 & 1.55 \\ \hline 
10 & 72.63 & 40.95 & 42.5 & 1.55 \\ \hline 
11 & 70.38 & 38.52 & 39.83 & 1.31 \\ \hline 
12 & 75.95 & 30.48 & 34.87 & 4.39 \\ \hline 
13 & 69.92 & 28.18 & 30.42 & 2.24 \\ \hline 
14 & 94.79 & 44.9 & 49.08 & 4.18 \\ \hline 
15 & 86.77 & 45.46 & 49.62 & 4.16 \\ \hline 
16 & 101.08 & 53.0 & 57.0 & 4.0 \\ \hline 
17 & 95.05 & 49.82 & 52.99 & 3.17
\\ \hline

\end{tabular}
\caption{Result of DIST strategy on Gaussian distributed customers with Deterministic travel time distribution assumption.}
}
\end{table}

\begin{table}
\parbox{.45\linewidth}{
\centering
\begin{tabular}{p{5mm}p{14mm}p{10mm}p{15mm}p{11mm}}
            \hline
            Id & Revenue & Served & Accepted & Missed \\ \hline
            
0 & 29.06 & 11.59 & 13.21 & 1.62 \\ \hline 
1 & 26.73 & 10.88 & 11.96 & 1.08 \\ \hline 
2 & 36.28 & 17.64 & 18.56 & 0.92 \\ \hline 
3 & 34.15 & 16.67 & 17.3 & 0.63 \\ \hline 
4 & 37.66 & 18.24 & 18.53 & 0.29 \\ \hline 
5 & 36.63 & 18.34 & 18.64 & 0.3 \\ \hline 
6 & 39.78 & 12.78 & 15.73 & 2.95 \\ \hline 
7 & 34.6 & 12.92 & 15.72 & 2.8 \\ \hline 
8 & 56.84 & 24.48 & 31.2 & 6.72 \\ \hline 
9 & 52.48 & 20.0 & 21.91 & 1.91 \\ \hline 
10 & 66.65 & 28.48 & 31.25 & 2.77 \\ \hline 
11 & 61.37 & 26.52 & 28.28 & 1.76 \\ \hline 
12 & 44.38 & 13.26 & 16.69 & 3.43 \\ \hline 
13 & 39.42 & 11.55 & 13.01 & 1.46 \\ \hline 
14 & 66.86 & 26.94 & 38.47 & 11.53 \\ \hline 
15 & 60.02 & 21.87 & 25.02 & 3.15 \\ \hline 
16 & 82.34 & 31.66 & 38.18 & 6.52 \\ \hline 
17 & 73.8 & 28.14 & 31.32 & 3.18
\\ \hline

\end{tabular}
\caption{Result of DIST strategy on Cluster distributed customers with Deterministic travel time distribution assumption.}
}
\hfill
\parbox{.45\linewidth}{
\centering
\begin{tabular}{p{5mm}p{14mm}p{10mm}p{15mm}p{11mm}}
            \hline
            Id & Revenue & Served & Accepted & Missed \\ \hline
            
0 & 27.35 & 11.27 & 13.32 & 2.05 \\ \hline 
1 & 24.58 & 10.63 & 11.96 & 1.33 \\ \hline 
2 & 34.15 & 17.76 & 19.18 & 1.42 \\ \hline 
3 & 32.21 & 17.36 & 18.51 & 1.15 \\ \hline 
4 & 36.35 & 19.09 & 19.68 & 0.59 \\ \hline 
5 & 35.39 & 18.69 & 19.22 & 0.53 \\ \hline 
6 & 38.12 & 13.0 & 17.07 & 4.07 \\ \hline 
7 & 34.02 & 11.17 & 12.32 & 1.15 \\ \hline 
8 & 53.09 & 23.51 & 30.84 & 7.33 \\ \hline 
9 & 49.39 & 19.15 & 20.78 & 1.63 \\ \hline 
10 & 62.13 & 30.39 & 35.9 & 5.51 \\ \hline 
11 & 58.91 & 26.1 & 28.07 & 1.97 \\ \hline 
12 & 42.46 & 13.62 & 20.38 & 6.76 \\ \hline 
13 & 38.8 & 12.4 & 14.34 & 1.94 \\ \hline 
14 & 65.13 & 24.7 & 32.8 & 8.1 \\ \hline 
15 & 55.38 & 20.87 & 23.99 & 3.12 \\ \hline 
16 & 78.15 & 32.17 & 40.47 & 8.3 \\ \hline 
17 & 70.39 & 29.24 & 33.05 & 3.81
\\ \hline

\end{tabular}
\caption{Result of DIST strategy on Uniform distributed customers with Deterministic travel time distribution assumption.}
}
\end{table}

\begin{table}
\parbox{.45\linewidth}{
\centering
\begin{tabular}{p{5mm}p{14mm}p{10mm}p{15mm}p{11mm}}
            \hline
            Id & Revenue & Served & Accepted & Missed \\ \hline
            
0 & 35.54 & 19.52 & 20.1 & 0.58 \\ \hline 
1 & 34.0 & 20.07 & 20.73 & 0.66 \\ \hline 
2 & 36.91 & 21.09 & 21.29 & 0.2 \\ \hline 
3 & 36.34 & 21.25 & 21.45 & 0.2 \\ \hline 
4 & 36.73 & 20.46 & 20.7 & 0.24 \\ \hline 
5 & 36.49 & 20.61 & 20.84 & 0.23 \\ \hline 
6 & 58.44 & 24.23 & 26.18 & 1.95 \\ \hline 
7 & 55.73 & 22.84 & 23.97 & 1.13 \\ \hline 
8 & 71.97 & 38.67 & 39.66 & 0.99 \\ \hline 
9 & 70.59 & 38.92 & 39.88 & 0.96 \\ \hline 
10 & 73.02 & 40.75 & 42.17 & 1.42 \\ \hline 
11 & 70.66 & 38.19 & 39.4 & 1.21 \\ \hline 
12 & 66.62 & 26.11 & 32.11 & 6.0 \\ \hline 
13 & 59.34 & 24.19 & 27.53 & 3.34 \\ \hline 
14 & 96.67 & 47.71 & 51.19 & 3.48 \\ \hline 
15 & 94.94 & 47.46 & 49.66 & 2.2 \\ \hline 
16 & 103.65 & 52.34 & 55.15 & 2.81 \\ \hline 
17 & 97.56 & 55.92 & 59.11 & 3.19
\\ \hline

\end{tabular}
\caption{Result of DIST strategy on Gaussian distributed customers with Gaussian travel time distribution assumption.}
}
\hfill
\parbox{.45\linewidth}{
\centering
\begin{tabular}{p{5mm}p{14mm}p{10mm}p{15mm}p{11mm}}
            \hline
            Id & Revenue & Served & Accepted & Missed \\ \hline
            
0 & 27.95 & 11.16 & 13.48 & 2.32 \\ \hline 
1 & 25.21 & 9.46 & 10.07 & 0.61 \\ \hline 
2 & 35.64 & 19.13 & 20.18 & 1.05 \\ \hline 
3 & 34.79 & 19.15 & 20.03 & 0.88 \\ \hline 
4 & 37.77 & 18.29 & 18.6 & 0.31 \\ \hline 
5 & 37.23 & 18.28 & 18.53 & 0.25 \\ \hline 
6 & 37.86 & 11.79 & 14.32 & 2.53 \\ \hline 
7 & 33.13 & 11.57 & 13.9 & 2.33 \\ \hline 
8 & 58.22 & 22.93 & 26.74 & 3.81 \\ \hline 
9 & 55.97 & 22.19 & 23.9 & 1.71 \\ \hline 
10 & 66.7 & 26.8 & 28.52 & 1.72 \\ \hline 
11 & 64.77 & 27.39 & 28.71 & 1.32 \\ \hline 
12 & 42.5 & 12.98 & 18.3 & 5.32 \\ \hline 
13 & 37.24 & 10.33 & 11.28 & 0.95 \\ \hline 
14 & 66.88 & 27.4 & 37.84 & 10.44 \\ \hline 
15 & 67.3 & 26.64 & 30.62 & 3.98 \\ \hline 
16 & 80.41 & 32.79 & 38.5 & 5.71 \\ \hline 
17 & 82.45 & 34.04 & 37.35 & 3.31
\\ \hline

\end{tabular}
\caption{Result of DIST strategy on Cluster distributed customers with Gaussian travel time distribution assumption.}
}
\end{table}

\begin{table}
\parbox{.45\linewidth}{
\centering
\begin{tabular}{p{5mm}p{14mm}p{10mm}p{15mm}p{11mm}}
            \hline
            Id & Revenue & Served & Accepted & Missed \\ \hline
            
0 & 26.15 & 10.65 & 12.72 & 2.07 \\ \hline 
1 & 23.59 & 10.16 & 11.58 & 1.42 \\ \hline 
2 & 34.51 & 17.41 & 18.48 & 1.07 \\ \hline 
3 & 33.31 & 17.29 & 18.1 & 0.81 \\ \hline 
4 & 36.36 & 19.13 & 19.7 & 0.57 \\ \hline 
5 & 35.39 & 19.18 & 19.71 & 0.53 \\ \hline 
6 & 36.12 & 12.79 & 17.48 & 4.69 \\ \hline 
7 & 32.07 & 11.59 & 13.73 & 2.14 \\ \hline 
8 & 53.27 & 21.25 & 24.47 & 3.22 \\ \hline 
9 & 50.0 & 23.36 & 26.75 & 3.39 \\ \hline 
10 & 63.48 & 27.77 & 30.58 & 2.81 \\ \hline 
11 & 60.99 & 28.37 & 30.51 & 2.14 \\ \hline 
12 & 40.39 & 12.96 & 18.36 & 5.4 \\ \hline 
13 & 36.64 & 10.58 & 11.52 & 0.94 \\ \hline 
14 & 64.49 & 23.0 & 27.43 & 4.43 \\ \hline 
15 & 63.33 & 24.18 & 27.14 & 2.96 \\ \hline 
16 & 76.44 & 31.96 & 39.5 & 7.54 \\ \hline 
17 & 75.46 & 34.76 & 39.58 & 4.82
\\ \hline

\end{tabular}
\caption{Result of DIST strategy on Uniform distributed customers with Gaussian travel time distribution assumption.}
}
\hfill
\parbox{.45\linewidth}{
\centering
\begin{tabular}{p{5mm}p{14mm}p{10mm}p{15mm}p{11mm}}
            \hline
            Id & Revenue & Served & Accepted & Missed \\ \hline
            
0 & 36.62 & 20.17 & 20.24 & 0.07 \\ \hline 
1 & 36.59 & 19.81 & 19.89 & 0.08 \\ \hline 
2 & 36.8 & 20.1 & 20.28 & 0.18 \\ \hline 
3 & 36.38 & 20.29 & 20.5 & 0.21 \\ \hline 
4 & 36.88 & 20.54 & 20.78 & 0.24 \\ \hline 
5 & 36.28 & 20.47 & 20.69 & 0.22 \\ \hline 
6 & 66.22 & 32.71 & 34.94 & 2.23 \\ \hline 
7 & 63.38 & 29.24 & 30.32 & 1.08 \\ \hline 
8 & 71.95 & 39.33 & 40.56 & 1.23 \\ \hline 
9 & 69.36 & 40.12 & 41.45 & 1.33 \\ \hline 
10 & 72.9 & 38.5 & 39.77 & 1.27 \\ \hline 
11 & 70.1 & 37.52 & 38.73 & 1.21 \\ \hline 
12 & 81.5 & 34.49 & 38.15 & 3.66 \\ \hline 
13 & 75.98 & 31.92 & 34.13 & 2.21 \\ \hline 
14 & 96.57 & 44.64 & 47.35 & 2.71 \\ \hline 
15 & 93.15 & 48.07 & 50.87 & 2.8 \\ \hline 
16 & 102.54 & 51.92 & 54.84 & 2.92 \\ \hline 
17 & 98.06 & 52.48 & 55.24 & 2.76
\\ \hline

\end{tabular}
\caption{Result of DIST strategy on Gaussian distributed customers with Misspecified travel time distribution assumption.}
}
\end{table}

\begin{table}
\parbox{.45\linewidth}{
\centering
\begin{tabular}{p{5mm}p{14mm}p{10mm}p{15mm}p{11mm}}
            \hline
            Id & Revenue & Served & Accepted & Missed \\ \hline
            
0 & 30.07 & 12.14 & 13.27 & 1.13 \\ \hline 
1 & 27.44 & 10.41 & 11.07 & 0.66 \\ \hline 
2 & 35.2 & 16.35 & 16.87 & 0.52 \\ \hline 
3 & 34.33 & 15.15 & 15.41 & 0.26 \\ \hline 
4 & 36.21 & 19.38 & 19.8 & 0.42 \\ \hline 
5 & 35.71 & 18.48 & 18.78 & 0.3 \\ \hline 
6 & 41.49 & 13.7 & 15.4 & 1.7 \\ \hline 
7 & 38.24 & 13.24 & 15.11 & 1.87 \\ \hline 
8 & 56.32 & 20.59 & 22.84 & 2.25 \\ \hline 
9 & 54.53 & 24.83 & 28.11 & 3.28 \\ \hline 
10 & 64.1 & 24.92 & 26.25 & 1.33 \\ \hline 
11 & 62.88 & 26.12 & 27.39 & 1.27 \\ \hline 
12 & 50.5 & 17.55 & 23.72 & 6.17 \\ \hline 
13 & 31.31 & 15.28 & 23.06 & 7.78 \\ \hline 
14 & 66.35 & 23.34 & 26.82 & 3.48 \\ \hline 
15 & 66.21 & 26.55 & 30.4 & 3.85 \\ \hline 
16 & 79.58 & 31.18 & 37.58 & 6.4 \\ \hline 
17 & 77.78 & 34.98 & 39.54 & 4.56
\\ \hline

\end{tabular}
\caption{Result of DIST strategy on Cluster distributed customers with Misspecified travel time distribution assumption.}
}
\hfill
\parbox{.45\linewidth}{
\centering
\begin{tabular}{p{5mm}p{14mm}p{10mm}p{15mm}p{11mm}}
            \hline
            Id & Revenue & Served & Accepted & Missed \\ \hline
            
0 & 28.37 & 11.71 & 13.13 & 1.42 \\ \hline 
1 & 26.11 & 11.12 & 12.15 & 1.03 \\ \hline 
2 & 34.33 & 16.48 & 17.27 & 0.79 \\ \hline 
3 & 33.05 & 16.75 & 17.53 & 0.78 \\ \hline 
4 & 35.56 & 19.99 & 20.65 & 0.66 \\ \hline 
5 & 34.4 & 20.01 & 20.62 & 0.61 \\ \hline 
6 & 41.21 & 14.94 & 19.16 & 4.22 \\ \hline 
7 & 36.66 & 13.1 & 14.8 & 1.7 \\ \hline 
8 & 51.94 & 22.2 & 28.12 & 5.92 \\ \hline 
9 & 51.91 & 21.54 & 23.62 & 2.08 \\ \hline 
10 & 62.45 & 27.49 & 30.64 & 3.15 \\ \hline 
11 & 58.99 & 29.66 & 32.66 & 3.0 \\ \hline 
12 & 46.22 & 16.33 & 24.73 & 8.4 \\ \hline 
13 & 26.67 & 13.88 & 21.65 & 7.77 \\ \hline 
14 & 61.25 & 23.53 & 32.32 & 8.79 \\ \hline 
15 & 62.04 & 25.75 & 30.19 & 4.44 \\ \hline 
16 & 75.62 & 30.04 & 36.65 & 6.61 \\ \hline 
17 & 69.6 & 35.62 & 42.53 & 6.91
\\ \hline

\end{tabular}
\caption{Result of DIST strategy on Uniform distributed customers with Misspecified travel time distribution assumption.}
}
\end{table}

\begin{table}
\parbox{.45\linewidth}{
\centering
\begin{tabular}{p{5mm}p{14mm}p{10mm}p{15mm}p{11mm}}
            \hline
            Id & Revenue & Served & Accepted & Missed \\ \hline
            
0 & 34.0 & 21.88 & 22.38 & 0.5 \\ \hline 
1 & 32.96 & 15.4 & 15.46 & 0.06 \\ \hline 
2 & 36.13 & 21.58 & 21.81 & 0.23 \\ \hline 
3 & 35.72 & 21.55 & 21.76 & 0.21 \\ \hline 
4 & 36.05 & 18.84 & 19.05 & 0.21 \\ \hline 
5 & 35.62 & 18.84 & 19.06 & 0.22 \\ \hline 
6 & 55.94 & 30.18 & 34.89 & 4.71 \\ \hline 
7 & 55.42 & 24.32 & 25.3 & 0.98 \\ \hline 
8 & 70.26 & 36.45 & 37.72 & 1.27 \\ \hline 
9 & 67.37 & 33.94 & 35.0 & 1.06 \\ \hline 
10 & 71.41 & 36.55 & 37.71 & 1.16 \\ \hline 
11 & 69.08 & 36.52 & 37.7 & 1.18 \\ \hline 
12 & 69.42 & 29.46 & 33.09 & 3.63 \\ \hline 
13 & 65.04 & 25.75 & 27.25 & 1.5 \\ \hline 
14 & 90.82 & 42.35 & 46.32 & 3.97 \\ \hline 
15 & 85.92 & 40.31 & 42.76 & 2.45 \\ \hline 
16 & 99.26 & 50.38 & 54.0 & 3.62 \\ \hline 
17 & 93.18 & 47.43 & 50.21 & 2.78
\\ \hline 

\end{tabular} \caption{Result of OPP strategy on Gaussian distributed customers with Deterministic travel time distribution assumption.}
}
\hfill
\parbox{.45\linewidth}{
\centering
\begin{tabular}{p{5mm}p{14mm}p{10mm}p{15mm}p{11mm}}
            \hline
            Id & Revenue & Served & Accepted & Missed \\ \hline
            
0 & 26.81 & 12.17 & 13.59 & 1.42 \\ \hline 
1 & 23.61 & 11.68 & 12.91 & 1.23 \\ \hline 
2 & 33.95 & 16.87 & 17.5 & 0.63 \\ \hline 
3 & 32.05 & 14.58 & 14.9 & 0.32 \\ \hline 
4 & 35.37 & 21.1 & 21.56 & 0.46 \\ \hline 
5 & 34.59 & 21.13 & 21.58 & 0.45 \\ \hline 
6 & 35.99 & 16.19 & 22.51 & 6.32 \\ \hline 
7 & 32.22 & 11.91 & 13.28 & 1.37 \\ \hline 
8 & 53.89 & 22.06 & 25.09 & 3.03 \\ \hline 
9 & 48.75 & 20.85 & 22.96 & 2.11 \\ \hline 
10 & 62.74 & 29.19 & 32.16 & 2.97 \\ \hline 
11 & 58.79 & 27.77 & 29.74 & 1.97 \\ \hline 
12 & 43.41 & 15.39 & 18.97 & 3.58 \\ \hline 
13 & 36.18 & 14.01 & 16.68 & 2.67 \\ \hline 
14 & 65.53 & 25.69 & 32.65 & 6.96 \\ \hline 
15 & 55.1 & 24.74 & 29.83 & 5.09 \\ \hline 
16 & 77.24 & 35.6 & 46.11 & 10.51 \\ \hline 
17 & 68.95 & 27.6 & 29.98 & 2.38
\\ \hline 

\end{tabular} \caption{Result of OPP strategy on Cluster distributed customers with Deterministic travel time distribution assumption.}
}
\end{table}

\begin{table}
\parbox{.45\linewidth}{
\centering
\begin{tabular}{p{5mm}p{14mm}p{10mm}p{15mm}p{11mm}}
            \hline
            Id & Revenue & Served & Accepted & Missed \\ \hline
            
0 & 24.77 & 11.36 & 13.48 & 2.12 \\ \hline 
1 & 21.32 & 10.69 & 12.15 & 1.46 \\ \hline 
2 & 33.19 & 16.51 & 17.57 & 1.06 \\ \hline 
3 & 31.3 & 14.55 & 15.09 & 0.54 \\ \hline 
4 & 35.28 & 18.4 & 18.94 & 0.54 \\ \hline 
5 & 33.58 & 20.84 & 21.56 & 0.72 \\ \hline 
6 & 33.26 & 13.05 & 17.37 & 4.32 \\ \hline 
7 & 27.17 & 12.41 & 15.34 & 2.93 \\ \hline 
8 & 50.85 & 22.75 & 29.55 & 6.8 \\ \hline 
9 & 45.45 & 18.82 & 20.89 & 2.07 \\ \hline 
10 & 60.11 & 28.0 & 32.26 & 4.26 \\ \hline 
11 & 55.0 & 26.94 & 29.73 & 2.79 \\ \hline 
12 & 39.0 & 13.9 & 18.26 & 4.36 \\ \hline 
13 & 33.89 & 11.95 & 13.87 & 1.92 \\ \hline 
14 & 61.67 & 25.22 & 35.35 & 10.13 \\ \hline 
15 & 52.32 & 19.77 & 22.55 & 2.78 \\ \hline 
16 & 73.73 & 30.98 & 39.36 & 8.38 \\ \hline 
17 & 62.66 & 30.18 & 36.08 & 5.9
\\ \hline 

\end{tabular} \caption{Result of OPP strategy on Uniform distributed customers with Deterministic travel time distribution assumption.}
}
\hfill
\parbox{.45\linewidth}{
\centering
\begin{tabular}{p{5mm}p{14mm}p{10mm}p{15mm}p{11mm}}
            \hline
            Id & Revenue & Served & Accepted & Missed \\ \hline
            
0 & 30.15 & 17.77 & 17.82 & 0.05 \\ \hline 
1 & 29.91 & 19.03 & 19.08 & 0.05 \\ \hline 
2 & 36.01 & 21.5 & 21.71 & 0.21 \\ \hline 
3 & 35.56 & 21.49 & 21.74 & 0.25 \\ \hline 
4 & 36.16 & 21.53 & 21.8 & 0.27 \\ \hline 
5 & 35.64 & 21.52 & 21.79 & 0.27 \\ \hline 
6 & 53.19 & 27.9 & 28.93 & 1.03 \\ \hline 
7 & 51.58 & 27.91 & 28.64 & 0.73 \\ \hline 
8 & 69.85 & 36.09 & 37.03 & 0.94 \\ \hline 
9 & 68.13 & 38.6 & 39.64 & 1.04 \\ \hline 
10 & 71.45 & 41.58 & 43.08 & 1.5 \\ \hline 
11 & 69.01 & 36.52 & 37.65 & 1.13 \\ \hline 
12 & 66.32 & 33.73 & 36.89 & 3.16 \\ \hline 
13 & 67.0 & 31.1 & 31.83 & 0.73 \\ \hline 
14 & 93.45 & 48.7 & 52.43 & 3.73 \\ \hline 
15 & 91.1 & 49.73 & 52.24 & 2.51 \\ \hline 
16 & 101.23 & 51.22 & 53.89 & 2.67 \\ \hline 
17 & 96.26 & 51.25 & 53.91 & 2.66
\\ \hline 

\end{tabular} \caption{Result of OPP strategy on Gaussian distributed customers with Gaussian travel time distribution assumption.}
}
\end{table}

\begin{table}
\parbox{.45\linewidth}{
\centering
\begin{tabular}{p{5mm}p{14mm}p{10mm}p{15mm}p{11mm}}
            \hline
            Id & Revenue & Served & Accepted & Missed \\ \hline
            
0 & 26.81 & 12.68 & 13.88 & 1.2 \\ \hline 
1 & 24.94 & 11.64 & 12.2 & 0.56 \\ \hline 
2 & 34.39 & 18.09 & 18.72 & 0.63 \\ \hline 
3 & 33.61 & 17.01 & 17.39 & 0.38 \\ \hline 
4 & 35.45 & 18.51 & 18.79 & 0.28 \\ \hline 
5 & 34.99 & 18.55 & 18.83 & 0.28 \\ \hline 
6 & 39.19 & 15.7 & 18.93 & 3.23 \\ \hline 
7 & 36.7 & 16.1 & 17.85 & 1.75 \\ \hline 
8 & 55.46 & 22.78 & 24.91 & 2.13 \\ \hline 
9 & 54.26 & 26.59 & 29.34 & 2.75 \\ \hline 
10 & 63.34 & 27.91 & 29.7 & 1.79 \\ \hline 
11 & 62.02 & 30.34 & 32.18 & 1.84 \\ \hline 
12 & 48.33 & 17.64 & 20.31 & 2.67 \\ \hline 
13 & 45.46 & 23.3 & 26.08 & 2.78 \\ \hline 
14 & 67.44 & 27.14 & 32.56 & 5.42 \\ \hline 
15 & 69.45 & 29.42 & 32.23 & 2.81 \\ \hline 
16 & 78.48 & 33.72 & 39.13 & 5.41 \\ \hline 
17 & 78.72 & 33.48 & 35.82 & 2.34
\\ \hline 

\end{tabular} \caption{Result of OPP strategy on Cluster distributed customers with Gaussian travel time distribution assumption.}
}
\hfill
\parbox{.45\linewidth}{
\centering
\begin{tabular}{p{5mm}p{14mm}p{10mm}p{15mm}p{11mm}}
            \hline
            Id & Revenue & Served & Accepted & Missed \\ \hline
            
0 & 23.58 & 12.61 & 15.71 & 3.1 \\ \hline 
1 & 22.4 & 10.22 & 10.99 & 0.77 \\ \hline 
2 & 33.11 & 15.58 & 16.25 & 0.67 \\ \hline 
3 & 32.32 & 15.69 & 16.21 & 0.52 \\ \hline 
4 & 35.37 & 19.65 & 20.23 & 0.58 \\ \hline 
5 & 34.19 & 17.17 & 17.53 & 0.36 \\ \hline 
6 & 36.01 & 14.56 & 19.13 & 4.57 \\ \hline 
7 & 32.73 & 15.07 & 17.31 & 2.24 \\ \hline 
8 & 51.48 & 20.4 & 22.71 & 2.31 \\ \hline 
9 & 50.71 & 22.51 & 24.65 & 2.14 \\ \hline 
10 & 60.92 & 28.5 & 32.05 & 3.55 \\ \hline 
11 & 57.99 & 29.39 & 32.07 & 2.68 \\ \hline 
12 & 42.09 & 16.76 & 23.06 & 6.3 \\ \hline 
13 & 39.21 & 19.62 & 22.5 & 2.88 \\ \hline 
14 & 61.54 & 22.86 & 27.27 & 4.41 \\ \hline 
15 & 62.73 & 27.93 & 31.95 & 4.02 \\ \hline 
16 & 73.14 & 30.11 & 36.14 & 6.03 \\ \hline 
17 & 73.1 & 34.64 & 39.21 & 4.57
\\ \hline 

\end{tabular} \caption{Result of OPP strategy on Uniform distributed customers with Gaussian travel time distribution assumption.}
}
\end{table}

\begin{table}
\parbox{.45\linewidth}{
\centering
\begin{tabular}{p{5mm}p{14mm}p{10mm}p{15mm}p{11mm}}
            \hline
            Id & Revenue & Served & Accepted & Missed \\ \hline
            
0 & 34.26 & 21.92 & 22.11 & 0.19 \\ \hline 
1 & 34.45 & 20.59 & 20.69 & 0.1 \\ \hline 
2 & 36.15 & 21.55 & 21.8 & 0.25 \\ \hline 
3 & 35.69 & 21.57 & 21.8 & 0.23 \\ \hline 
4 & 36.13 & 21.51 & 21.78 & 0.27 \\ \hline 
5 & 35.57 & 21.49 & 21.77 & 0.28 \\ \hline 
6 & 61.3 & 30.62 & 33.07 & 2.45 \\ \hline 
7 & 60.4 & 29.59 & 30.71 & 1.12 \\ \hline 
8 & 70.88 & 36.63 & 37.65 & 1.02 \\ \hline 
9 & 69.04 & 36.68 & 37.68 & 1.0 \\ \hline 
10 & 71.55 & 41.66 & 43.12 & 1.46 \\ \hline 
11 & 69.06 & 36.51 & 37.64 & 1.13 \\ \hline 
12 & 75.91 & 33.15 & 36.73 & 3.58 \\ \hline 
13 & 73.03 & 33.48 & 35.56 & 2.08 \\ \hline 
14 & 93.86 & 43.74 & 46.46 & 2.72 \\ \hline 
15 & 91.08 & 41.4 & 42.86 & 1.46 \\ \hline 
16 & 100.76 & 51.06 & 53.98 & 2.92 \\ \hline 
17 & 95.88 & 51.25 & 53.98 & 2.73
\\ \hline 

\end{tabular} \caption{Result of OPP strategy on Gaussian distributed customers with Misspecified travel time distribution assumption.}
}
\hfill
\parbox{.45\linewidth}{
\centering
\begin{tabular}{p{5mm}p{14mm}p{10mm}p{15mm}p{11mm}}
            \hline
            Id & Revenue & Served & Accepted & Missed \\ \hline
            
0 & 27.74 & 13.2 & 15.11 & 1.91 \\ \hline 
1 & 25.27 & 11.42 & 12.09 & 0.67 \\ \hline 
2 & 34.33 & 18.1 & 18.84 & 0.74 \\ \hline 
3 & 32.43 & 14.71 & 14.96 & 0.25 \\ \hline 
4 & 35.67 & 19.84 & 20.17 & 0.33 \\ \hline 
5 & 34.6 & 17.28 & 17.46 & 0.18 \\ \hline 
6 & 40.86 & 16.5 & 21.05 & 4.55 \\ \hline 
7 & 33.91 & 17.74 & 21.16 & 3.42 \\ \hline 
8 & 54.36 & 22.37 & 25.03 & 2.66 \\ \hline 
9 & 53.01 & 23.1 & 25.03 & 1.93 \\ \hline 
10 & 63.16 & 29.44 & 32.16 & 2.72 \\ \hline 
11 & 60.95 & 28.17 & 29.7 & 1.53 \\ \hline 
12 & 49.17 & 18.68 & 24.36 & 5.68 \\ \hline 
13 & 45.81 & 16.37 & 17.78 & 1.41 \\ \hline 
14 & 64.56 & 24.94 & 30.03 & 5.09 \\ \hline 
15 & 64.31 & 28.42 & 32.67 & 4.25 \\ \hline 
16 & 76.15 & 31.17 & 35.93 & 4.76 \\ \hline 
17 & 75.53 & 32.77 & 35.84 & 3.07
\\ \hline 

\end{tabular} \caption{Result of OPP strategy on Cluster distributed customers with Misspecified travel time distribution assumption.}
}
\end{table}

\begin{table}
\parbox{.45\linewidth}{
\centering
\begin{tabular}{p{5mm}p{14mm}p{10mm}p{15mm}p{11mm}}
            \hline
            Id & Revenue & Served & Accepted & Missed \\ \hline
            
0 & 26.08 & 11.84 & 13.97 & 2.13 \\ \hline 
1 & 23.94 & 11.14 & 12.18 & 1.04 \\ \hline 
2 & 33.29 & 17.63 & 18.93 & 1.3 \\ \hline 
3 & 32.49 & 18.06 & 18.93 & 0.87 \\ \hline 
4 & 35.43 & 18.5 & 19.0 & 0.5 \\ \hline 
5 & 34.14 & 20.95 & 21.5 & 0.55 \\ \hline 
6 & 36.24 & 14.65 & 20.73 & 6.08 \\ \hline 
7 & 33.05 & 13.82 & 15.99 & 2.17 \\ \hline 
8 & 50.07 & 21.74 & 27.26 & 5.52 \\ \hline 
9 & 48.75 & 22.14 & 24.97 & 2.83 \\ \hline 
10 & 60.0 & 28.12 & 32.23 & 4.11 \\ \hline 
11 & 57.67 & 27.46 & 29.72 & 2.26 \\ \hline 
12 & 42.92 & 15.12 & 19.97 & 4.85 \\ \hline 
13 & 35.96 & 17.65 & 21.92 & 4.27 \\ \hline 
14 & 58.12 & 23.3 & 33.01 & 9.71 \\ \hline 
15 & 57.98 & 25.29 & 29.93 & 4.64 \\ \hline 
16 & 71.29 & 29.2 & 36.14 & 6.94 \\ \hline 
17 & 70.16 & 27.84 & 30.14 & 2.3
\\ \hline 

\end{tabular} \caption{Result of OPP strategy on Uniform distributed customers with Misspecified travel time distribution assumption.}
}
\hfill
\parbox{.45\linewidth}{
\centering
\begin{tabular}{p{5mm}p{14mm}p{10mm}p{15mm}p{11mm}}
            \hline
            Id & Revenue & Served & Accepted & Missed \\ \hline
            
0 & 36.77 & 18.61 & 19.01 & 0.4 \\ \hline 
1 & 35.22 & 17.7 & 18.0 & 0.3 \\ \hline 
2 & 38.9 & 19.88 & 20.08 & 0.2 \\ \hline 
3 & 38.56 & 19.92 & 20.13 & 0.21 \\ \hline 
4 & 38.99 & 19.91 & 20.12 & 0.21 \\ \hline 
5 & 38.56 & 19.9 & 20.12 & 0.22 \\ \hline 
6 & 54.54 & 26.47 & 32.38 & 5.91 \\ \hline 
7 & 49.34 & 23.93 & 26.12 & 2.19 \\ \hline 
8 & 74.96 & 37.81 & 39.92 & 2.11 \\ \hline 
9 & 71.34 & 37.98 & 39.93 & 1.95 \\ \hline 
10 & 77.22 & 38.54 & 39.93 & 1.39 \\ \hline 
11 & 74.4 & 38.52 & 39.92 & 1.4 \\ \hline 
12 & 59.37 & 28.81 & 43.29 & 14.48 \\ \hline 
13 & 50.69 & 27.12 & 32.97 & 5.85 \\ \hline 
14 & 91.39 & 45.93 & 57.18 & 11.25 \\ \hline 
15 & 72.96 & 46.94 & 57.06 & 10.12 \\ \hline 
16 & 105.74 & 52.2 & 57.28 & 5.08 \\ \hline 
17 & 96.72 & 52.44 & 57.28 & 4.84
\\ \hline 

\end{tabular} \caption{Result of OPT strategy on Gaussian distributed customers with Deterministic travel time distribution assumption.}
}
\end{table}

\begin{table}
\parbox{.45\linewidth}{
\centering
\begin{tabular}{p{5mm}p{14mm}p{10mm}p{15mm}p{11mm}}
            \hline
            Id & Revenue & Served & Accepted & Missed \\ \hline
            
0 & 25.43 & 12.36 & 15.78 & 3.42 \\ \hline 
1 & 21.73 & 10.82 & 12.3 & 1.48 \\ \hline 
2 & 36.17 & 18.47 & 19.97 & 1.5 \\ \hline 
3 & 33.41 & 18.51 & 19.93 & 1.42 \\ \hline 
4 & 38.35 & 19.55 & 19.99 & 0.44 \\ \hline 
5 & 37.52 & 19.55 & 19.98 & 0.43 \\ \hline 
6 & 31.21 & 16.06 & 27.95 & 11.89 \\ \hline 
7 & 20.65 & 14.22 & 19.49 & 5.27 \\ \hline 
8 & 50.57 & 25.96 & 39.1 & 13.14 \\ \hline 
9 & 28.34 & 26.61 & 38.07 & 11.46 \\ \hline 
10 & 64.31 & 32.33 & 39.64 & 7.31 \\ \hline 
11 & 51.9 & 32.92 & 39.64 & 6.72 \\ \hline 
12 & 34.42 & 18.78 & 35.82 & 17.04 \\ \hline 
13 & 19.26 & 15.81 & 23.56 & 7.75 \\ \hline 
14 & 55.17 & 29.14 & 53.54 & 24.4 \\ \hline 
15 & 17.35 & 30.41 & 49.79 & 19.38 \\ \hline 
16 & 74.28 & 37.39 & 56.58 & 19.19 \\ \hline 
17 & 41.08 & 38.42 & 55.67 & 17.25
\\ \hline 

\end{tabular} \caption{Result of OPT strategy on Cluster distributed customers with Deterministic travel time distribution assumption.}
}
\hfill
\parbox{.45\linewidth}{
\centering
\begin{tabular}{p{5mm}p{14mm}p{10mm}p{15mm}p{11mm}}
            \hline
            Id & Revenue & Served & Accepted & Missed \\ \hline
            
0 & 23.51 & 11.26 & 15.22 & 3.96 \\ \hline 
1 & 19.6 & 10.45 & 12.3 & 1.85 \\ \hline 
2 & 34.8 & 17.78 & 20.0 & 2.22 \\ \hline 
3 & 30.61 & 17.82 & 19.95 & 2.13 \\ \hline 
4 & 37.85 & 19.25 & 20.01 & 0.76 \\ \hline 
5 & 36.27 & 19.23 & 20.01 & 0.78 \\ \hline 
6 & 27.22 & 13.85 & 25.62 & 11.77 \\ \hline 
7 & 20.99 & 12.81 & 17.14 & 4.33 \\ \hline 
8 & 45.2 & 23.46 & 38.49 & 15.03 \\ \hline 
9 & 20.66 & 24.16 & 36.78 & 12.62 \\ \hline 
10 & 59.84 & 30.18 & 39.66 & 9.48 \\ \hline 
11 & 44.08 & 30.98 & 39.54 & 8.56 \\ \hline 
12 & 30.47 & 15.95 & 30.91 & 14.96 \\ \hline 
13 & 16.67 & 13.96 & 20.88 & 6.92 \\ \hline 
14 & 47.08 & 25.22 & 51.18 & 25.96 \\ \hline 
15 & 9.4 & 26.3 & 45.33 & 19.03 \\ \hline 
16 & 66.39 & 33.83 & 56.03 & 22.2 \\ \hline 
17 & 28.76 & 34.94 & 54.22 & 19.28
\\ \hline 

\end{tabular} \caption{Result of OPT strategy on Uniform distributed customers with Deterministic travel time distribution assumption.}
}
\end{table}

\begin{table}
\parbox{.45\linewidth}{
\centering
\begin{tabular}{p{5mm}p{14mm}p{10mm}p{15mm}p{11mm}}
            \hline
            Id & Revenue & Served & Accepted & Missed \\ \hline
            
0 & 32.92 & 16.65 & 16.7 & 0.05 \\ \hline 
1 & 33.41 & 16.99 & 17.02 & 0.03 \\ \hline 
2 & 38.78 & 19.79 & 20.01 & 0.22 \\ \hline 
3 & 38.44 & 19.84 & 20.05 & 0.21 \\ \hline 
4 & 39.0 & 19.91 & 20.12 & 0.21 \\ \hline 
5 & 38.58 & 19.91 & 20.15 & 0.24 \\ \hline 
6 & 55.63 & 28.1 & 29.54 & 1.44 \\ \hline 
7 & 54.59 & 27.98 & 28.7 & 0.72 \\ \hline 
8 & 75.18 & 37.95 & 39.1 & 1.15 \\ \hline 
9 & 73.31 & 38.03 & 39.09 & 1.06 \\ \hline 
10 & 77.31 & 38.56 & 39.9 & 1.34 \\ \hline 
11 & 74.74 & 38.61 & 39.92 & 1.31 \\ \hline 
12 & 65.2 & 33.35 & 38.03 & 4.68 \\ \hline 
13 & 63.34 & 34.2 & 36.19 & 1.99 \\ \hline 
14 & 98.04 & 49.8 & 54.57 & 4.77 \\ \hline 
15 & 96.16 & 51.43 & 54.45 & 3.02 \\ \hline 
16 & 108.73 & 53.59 & 56.94 & 3.35 \\ \hline 
17 & 103.18 & 53.7 & 56.77 & 3.07
\\ \hline 

\end{tabular} \caption{Result of OPT strategy on Gaussian distributed customers with Gaussian travel time distribution assumption.}
}
\hfill
\parbox{.45\linewidth}{
\centering
\begin{tabular}{p{5mm}p{14mm}p{10mm}p{15mm}p{11mm}}
            \hline
            Id & Revenue & Served & Accepted & Missed \\ \hline
            
0 & 26.68 & 13.31 & 15.39 & 2.08 \\ \hline 
1 & 25.02 & 13.4 & 14.4 & 1.0 \\ \hline 
2 & 36.74 & 18.84 & 19.83 & 0.99 \\ \hline 
3 & 35.69 & 19.02 & 19.76 & 0.74 \\ \hline 
4 & 38.54 & 19.66 & 19.99 & 0.33 \\ \hline 
5 & 37.84 & 19.62 & 19.94 & 0.32 \\ \hline 
6 & 36.01 & 18.92 & 26.28 & 7.36 \\ \hline 
7 & 33.12 & 20.33 & 23.03 & 2.7 \\ \hline 
8 & 52.89 & 27.86 & 36.81 & 8.95 \\ \hline 
9 & 48.82 & 31.38 & 37.22 & 5.84 \\ \hline 
10 & 65.01 & 33.23 & 38.59 & 5.36 \\ \hline 
11 & 62.46 & 35.27 & 39.08 & 3.81 \\ \hline 
12 & 44.03 & 23.74 & 33.91 & 10.17 \\ \hline 
13 & 38.91 & 25.61 & 29.13 & 3.52 \\ \hline 
14 & 62.32 & 33.86 & 49.44 & 15.58 \\ \hline 
15 & 56.74 & 40.65 & 50.02 & 9.37 \\ \hline 
16 & 75.71 & 39.76 & 52.2 & 12.44 \\ \hline 
17 & 73.58 & 45.87 & 54.4 & 8.53
\\ \hline 

\end{tabular} \caption{Result of OPT strategy on Cluster distributed customers with Gaussian travel time distribution assumption.}
}
\end{table}

\begin{table}
\parbox{.45\linewidth}{
\centering
\begin{tabular}{p{5mm}p{14mm}p{10mm}p{15mm}p{11mm}}
            \hline
            Id & Revenue & Served & Accepted & Missed \\ \hline
            
0 & 23.77 & 11.81 & 14.57 & 2.76 \\ \hline 
1 & 22.32 & 11.86 & 12.99 & 1.13 \\ \hline 
2 & 35.17 & 18.08 & 19.8 & 1.72 \\ \hline 
3 & 33.34 & 18.41 & 19.7 & 1.29 \\ \hline 
4 & 37.95 & 19.27 & 19.98 & 0.71 \\ \hline 
5 & 36.87 & 19.37 & 20.0 & 0.63 \\ \hline 
6 & 32.62 & 16.85 & 24.87 & 8.02 \\ \hline 
7 & 28.97 & 17.58 & 20.41 & 2.83 \\ \hline 
8 & 46.99 & 24.89 & 36.03 & 11.14 \\ \hline 
9 & 39.45 & 28.7 & 36.21 & 7.51 \\ \hline 
10 & 60.44 & 31.04 & 38.1 & 7.06 \\ \hline 
11 & 55.21 & 33.42 & 38.83 & 5.41 \\ \hline 
12 & 37.45 & 19.48 & 30.32 & 10.84 \\ \hline 
13 & 33.54 & 21.11 & 24.56 & 3.45 \\ \hline 
14 & 52.2 & 28.42 & 47.76 & 19.34 \\ \hline 
15 & 41.69 & 35.73 & 47.27 & 11.54 \\ \hline 
16 & 66.22 & 34.79 & 50.79 & 16.0 \\ \hline 
17 & 60.41 & 42.35 & 53.42 & 11.07
\\ \hline 

\end{tabular} \caption{Result of OPT strategy on Uniform distributed customers with Gaussian travel time distribution assumption.}
}
\hfill
\parbox{.45\linewidth}{
\centering
\begin{tabular}{p{5mm}p{14mm}p{10mm}p{15mm}p{11mm}}
            \hline
            Id & Revenue & Served & Accepted & Missed \\ \hline
            
0 & 37.45 & 19.1 & 19.2 & 0.1 \\ \hline 
1 & 37.29 & 19.09 & 19.17 & 0.08 \\ \hline 
2 & 39.0 & 19.91 & 20.12 & 0.21 \\ \hline 
3 & 38.58 & 19.91 & 20.1 & 0.19 \\ \hline 
4 & 38.98 & 19.9 & 20.11 & 0.21 \\ \hline 
5 & 38.65 & 19.95 & 20.16 & 0.21 \\ \hline 
6 & 59.6 & 30.05 & 34.11 & 4.06 \\ \hline 
7 & 55.87 & 28.48 & 30.22 & 1.74 \\ \hline 
8 & 76.28 & 38.48 & 39.88 & 1.4 \\ \hline 
9 & 74.1 & 38.66 & 39.91 & 1.25 \\ \hline 
10 & 77.35 & 38.62 & 39.94 & 1.32 \\ \hline 
11 & 74.68 & 38.59 & 39.9 & 1.31 \\ \hline 
12 & 68.4 & 33.78 & 44.63 & 10.85 \\ \hline 
13 & 63.56 & 33.38 & 37.67 & 4.29 \\ \hline 
14 & 94.35 & 48.21 & 55.65 & 7.44 \\ \hline 
15 & 90.58 & 51.29 & 56.94 & 5.65 \\ \hline 
16 & 107.99 & 53.38 & 57.03 & 3.65 \\ \hline 
17 & 102.56 & 53.88 & 57.24 & 3.36
\\ \hline 

\end{tabular} \caption{Result of OPT strategy on Gaussian distributed customers with Misspecified travel time distribution assumption.}
}
\end{table}

\begin{table}
\parbox{.45\linewidth}{
\centering
\begin{tabular}{p{5mm}p{14mm}p{10mm}p{15mm}p{11mm}}
            \hline
            Id & Revenue & Served & Accepted & Missed \\ \hline
            
0 & 27.53 & 13.93 & 16.77 & 2.84 \\ \hline 
1 & 25.16 & 13.2 & 14.45 & 1.25 \\ \hline 
2 & 36.69 & 18.82 & 19.92 & 1.1 \\ \hline 
3 & 35.69 & 19.15 & 19.99 & 0.84 \\ \hline 
4 & 38.61 & 19.68 & 20.05 & 0.37 \\ \hline 
5 & 37.95 & 19.69 & 20.05 & 0.36 \\ \hline 
6 & 36.95 & 18.97 & 28.53 & 9.56 \\ \hline 
7 & 29.89 & 18.39 & 22.0 & 3.61 \\ \hline 
8 & 49.82 & 26.34 & 36.37 & 10.03 \\ \hline 
9 & 43.05 & 30.37 & 38.23 & 7.86 \\ \hline 
10 & 63.32 & 32.34 & 38.31 & 5.97 \\ \hline 
11 & 59.51 & 34.63 & 39.21 & 4.58 \\ \hline 
12 & 42.24 & 21.83 & 35.83 & 14.0 \\ \hline 
13 & 34.07 & 22.06 & 27.31 & 5.25 \\ \hline 
14 & 56.07 & 30.48 & 49.64 & 19.16 \\ \hline 
15 & 40.49 & 37.51 & 52.3 & 14.79 \\ \hline 
16 & 69.86 & 36.74 & 51.36 & 14.62 \\ \hline 
17 & 62.23 & 43.46 & 55.32 & 11.86
\\ \hline 

\end{tabular} \caption{Result of OPT strategy on Cluster distributed customers with Misspecified travel time distribution assumption.}
}
\hfill
\parbox{.45\linewidth}{
\centering
\begin{tabular}{p{5mm}p{14mm}p{10mm}p{15mm}p{11mm}}
            \hline
            Id & Revenue & Served & Accepted & Missed \\ \hline
            
0 & 25.18 & 12.44 & 15.86 & 3.42 \\ \hline 
1 & 22.18 & 12.05 & 13.46 & 1.41 \\ \hline 
2 & 34.92 & 17.97 & 19.81 & 1.84 \\ \hline 
3 & 32.95 & 18.39 & 19.86 & 1.47 \\ \hline 
4 & 37.97 & 19.29 & 19.99 & 0.7 \\ \hline 
5 & 36.85 & 19.35 & 19.95 & 0.6 \\ \hline 
6 & 32.33 & 16.56 & 26.52 & 9.96 \\ \hline 
7 & 26.67 & 16.32 & 20.09 & 3.77 \\ \hline 
8 & 44.19 & 23.51 & 35.22 & 11.71 \\ \hline 
9 & 31.36 & 27.2 & 37.41 & 10.21 \\ \hline 
10 & 58.51 & 30.09 & 37.61 & 7.52 \\ \hline 
11 & 51.11 & 32.55 & 39.14 & 6.59 \\ \hline 
12 & 36.2 & 18.19 & 30.54 & 12.35 \\ \hline 
13 & 30.38 & 18.74 & 23.46 & 4.72 \\ \hline 
14 & 47.26 & 25.91 & 48.3 & 22.39 \\ \hline 
15 & 24.31 & 32.11 & 49.17 & 17.06 \\ \hline 
16 & 62.05 & 32.9 & 50.09 & 17.19 \\ \hline 
17 & 47.11 & 39.64 & 54.61 & 14.97
\\ \hline 

\end{tabular} \caption{Result of OPT strategy on Uniform distributed customers with Misspecified travel time distribution assumption.}
}
\end{table}

\begin{table}
\parbox{.45\linewidth}{
\centering
\begin{tabular}{p{5mm}p{14mm}p{10mm}p{15mm}p{11mm}}
            \hline
            Id & Revenue & Served & Accepted & Missed \\ \hline
            
0 & 35.64 & 16.86 & 17.05 & 0.19 \\ \hline 
1 & 33.06 & 14.63 & 14.7 & 0.07 \\ \hline 
2 & 38.74 & 18.76 & 18.95 & 0.19 \\ \hline 
3 & 38.22 & 18.61 & 18.82 & 0.21 \\ \hline 
4 & 38.65 & 18.65 & 18.81 & 0.16 \\ \hline 
5 & 38.49 & 18.81 & 18.97 & 0.16 \\ \hline 
6 & 59.35 & 26.19 & 28.75 & 2.56 \\ \hline 
7 & 55.55 & 24.25 & 25.64 & 1.39 \\ \hline 
8 & 74.72 & 35.37 & 36.79 & 1.42 \\ \hline 
9 & 71.95 & 34.76 & 35.98 & 1.22 \\ \hline 
10 & 76.86 & 36.61 & 37.79 & 1.18 \\ \hline 
11 & 74.45 & 36.6 & 37.78 & 1.18 \\ \hline 
12 & 72.42 & 30.19 & 34.45 & 4.26 \\ \hline 
13 & 64.96 & 26.77 & 28.91 & 2.14 \\ \hline 
14 & 95.79 & 44.84 & 51.86 & 7.02 \\ \hline 
15 & 87.95 & 41.02 & 44.36 & 3.34 \\ \hline 
16 & 106.13 & 49.46 & 53.38 & 3.92 \\ \hline 
17 & 99.82 & 48.62 & 51.8 & 3.18
\\ \hline 

\end{tabular} \caption{Result of OPT+basis strategy on Gaussian distributed customers with Deterministic travel time distribution assumption.}
}
\hfill
\parbox{.45\linewidth}{
\centering
\begin{tabular}{p{5mm}p{14mm}p{10mm}p{15mm}p{11mm}}
            \hline
            Id & Revenue & Served & Accepted & Missed \\ \hline
            
0 & 27.51 & 12.02 & 14.06 & 2.04 \\ \hline 
1 & 22.38 & 10.08 & 11.17 & 1.09 \\ \hline 
2 & 36.36 & 18.25 & 19.56 & 1.31 \\ \hline 
3 & 34.9 & 16.54 & 17.26 & 0.72 \\ \hline 
4 & 38.23 & 19.44 & 19.91 & 0.47 \\ \hline 
5 & 37.58 & 19.57 & 19.95 & 0.38 \\ \hline 
6 & 37.1 & 11.35 & 12.48 & 1.13 \\ \hline 
7 & 32.43 & 12.36 & 14.27 & 1.91 \\ \hline 
8 & 57.86 & 22.13 & 25.76 & 3.63 \\ \hline 
9 & 52.08 & 19.51 & 21.21 & 1.7 \\ \hline 
10 & 66.71 & 30.67 & 35.29 & 4.62 \\ \hline 
11 & 62.5 & 26.54 & 28.4 & 1.86 \\ \hline 
12 & 43.27 & 15.29 & 20.43 & 5.14 \\ \hline 
13 & 36.43 & 13.3 & 15.58 & 2.28 \\ \hline 
14 & 68.39 & 26.86 & 36.74 & 9.88 \\ \hline 
15 & 57.32 & 24.75 & 30.76 & 6.01 \\ \hline 
16 & 81.99 & 32.66 & 40.61 & 7.95 \\ \hline 
17 & 73.61 & 27.59 & 30.4 & 2.81
\\ \hline 

\end{tabular} \caption{Result of OPT+basis strategy on Cluster distributed customers with Deterministic travel time distribution assumption.}
}
\end{table}

\begin{table}
\parbox{.45\linewidth}{
\centering
\begin{tabular}{p{5mm}p{14mm}p{10mm}p{15mm}p{11mm}}
            \hline
            Id & Revenue & Served & Accepted & Missed \\ \hline
            
0 & 26.69 & 11.0 & 12.62 & 1.62 \\ \hline 
1 & 22.83 & 9.06 & 9.78 & 0.72 \\ \hline 
2 & 34.82 & 17.68 & 19.85 & 2.17 \\ \hline 
3 & 33.2 & 16.63 & 17.8 & 1.17 \\ \hline 
4 & 37.91 & 19.07 & 19.76 & 0.69 \\ \hline 
5 & 36.49 & 18.73 & 19.4 & 0.67 \\ \hline 
6 & 37.95 & 13.28 & 15.86 & 2.58 \\ \hline 
7 & 32.11 & 11.92 & 13.63 & 1.71 \\ \hline 
8 & 54.01 & 23.54 & 31.02 & 7.48 \\ \hline 
9 & 49.72 & 20.36 & 22.68 & 2.32 \\ \hline 
10 & 63.41 & 29.65 & 35.85 & 6.2 \\ \hline 
11 & 59.58 & 25.23 & 27.11 & 1.88 \\ \hline 
12 & 41.59 & 14.62 & 18.93 & 4.31 \\ \hline 
13 & 37.46 & 12.98 & 14.92 & 1.94 \\ \hline 
14 & 66.26 & 24.59 & 31.79 & 7.2 \\ \hline 
15 & 58.26 & 21.49 & 24.35 & 2.86 \\ \hline 
16 & 79.26 & 31.5 & 38.99 & 7.49 \\ \hline 
17 & 71.28 & 27.35 & 30.04 & 2.69
\\ \hline 

\end{tabular} \caption{Result of OPT+basis strategy on Uniform distributed customers with Deterministic travel time distribution assumption.}
}
\hfill
\parbox{.45\linewidth}{
\centering
\begin{tabular}{p{5mm}p{14mm}p{10mm}p{15mm}p{11mm}}
            \hline
            Id & Revenue & Served & Accepted & Missed \\ \hline
            
0 & 32.68 & 16.43 & 16.47 & 0.04 \\ \hline 
1 & 32.28 & 15.43 & 15.45 & 0.02 \\ \hline 
2 & 38.92 & 19.86 & 20.08 & 0.22 \\ \hline 
3 & 38.58 & 19.87 & 20.07 & 0.2 \\ \hline 
4 & 39.02 & 19.84 & 20.04 & 0.2 \\ \hline 
5 & 38.63 & 19.94 & 20.15 & 0.21 \\ \hline 
6 & 56.65 & 27.1 & 27.91 & 0.81 \\ \hline 
7 & 55.52 & 27.33 & 27.96 & 0.63 \\ \hline 
8 & 75.6 & 37.57 & 38.58 & 1.01 \\ \hline 
9 & 73.46 & 38.15 & 39.18 & 1.03 \\ \hline 
10 & 77.26 & 38.57 & 39.88 & 1.31 \\ \hline 
11 & 74.87 & 38.6 & 39.88 & 1.28 \\ \hline 
12 & 65.66 & 33.49 & 37.91 & 4.42 \\ \hline 
13 & 69.61 & 32.69 & 33.65 & 0.96 \\ \hline 
14 & 100.37 & 47.54 & 50.71 & 3.17 \\ \hline 
15 & 97.92 & 49.64 & 52.24 & 2.6 \\ \hline 
16 & 108.97 & 51.8 & 54.69 & 2.89 \\ \hline 
17 & 103.5 & 53.22 & 56.22 & 3.0
\\ \hline 

\end{tabular} \caption{Result of OPT+basis strategy on Gaussian distributed customers with Gaussian travel time distribution assumption.}
}
\end{table}

\begin{table}
\parbox{.45\linewidth}{
\centering
\begin{tabular}{p{5mm}p{14mm}p{10mm}p{15mm}p{11mm}}
            \hline
            Id & Revenue & Served & Accepted & Missed \\ \hline
            
0 & 26.88 & 10.67 & 11.27 & 0.6 \\ \hline 
1 & 25.76 & 12.14 & 12.86 & 0.72 \\ \hline 
2 & 36.68 & 18.47 & 19.4 & 0.93 \\ \hline 
3 & 36.08 & 17.91 & 18.42 & 0.51 \\ \hline 
4 & 38.47 & 19.61 & 19.95 & 0.34 \\ \hline 
5 & 37.85 & 19.62 & 19.94 & 0.32 \\ \hline 
6 & 40.38 & 15.85 & 19.55 & 3.7 \\ \hline 
7 & 36.72 & 15.8 & 17.4 & 1.6 \\ \hline 
8 & 58.49 & 22.84 & 25.5 & 2.66 \\ \hline 
9 & 56.45 & 22.18 & 23.65 & 1.47 \\ \hline 
10 & 67.76 & 30.18 & 32.97 & 2.79 \\ \hline 
11 & 64.43 & 33.17 & 36.13 & 2.96 \\ \hline 
12 & 48.67 & 19.71 & 25.36 & 5.65 \\ \hline 
13 & 47.24 & 18.02 & 19.63 & 1.61 \\ \hline 
14 & 70.09 & 26.88 & 32.45 & 5.57 \\ \hline 
15 & 71.58 & 28.86 & 32.28 & 3.42 \\ \hline 
16 & 81.83 & 33.71 & 39.49 & 5.78 \\ \hline 
17 & 82.85 & 36.36 & 39.95 & 3.59
\\ \hline 

\end{tabular} \caption{Result of OPT+basis strategy on Cluster distributed customers with Gaussian travel time distribution assumption.}
}
\hfill
\parbox{.45\linewidth}{
\centering
\begin{tabular}{p{5mm}p{14mm}p{10mm}p{15mm}p{11mm}}
            \hline
            Id & Revenue & Served & Accepted & Missed \\ \hline
            
0 & 26.72 & 11.51 & 12.94 & 1.43 \\ \hline 
1 & 24.82 & 10.45 & 11.05 & 0.6 \\ \hline 
2 & 35.31 & 16.33 & 17.21 & 0.88 \\ \hline 
3 & 34.12 & 17.45 & 18.44 & 0.99 \\ \hline 
4 & 38.03 & 18.95 & 19.57 & 0.62 \\ \hline 
5 & 36.97 & 19.39 & 20.01 & 0.62 \\ \hline 
6 & 37.19 & 15.33 & 19.77 & 4.44 \\ \hline 
7 & 35.15 & 14.72 & 16.46 & 1.74 \\ \hline 
8 & 55.63 & 22.45 & 25.73 & 3.28 \\ \hline 
9 & 53.44 & 20.88 & 22.31 & 1.43 \\ \hline 
10 & 63.47 & 30.14 & 35.09 & 4.95 \\ \hline 
11 & 62.22 & 28.76 & 31.1 & 2.34 \\ \hline 
12 & 45.75 & 17.86 & 22.62 & 4.76 \\ \hline 
13 & 44.08 & 15.4 & 16.43 & 1.03 \\ \hline 
14 & 66.04 & 24.19 & 28.59 & 4.4 \\ \hline 
15 & 66.43 & 26.46 & 29.49 & 3.03 \\ \hline 
16 & 78.24 & 31.79 & 37.63 & 5.84 \\ \hline 
17 & 77.72 & 32.61 & 36.0 & 3.39
\\ \hline 

\end{tabular} \caption{Result of OPT+basis strategy on Uniform distributed customers with Gaussian travel time distribution assumption.}
}
\end{table}

\begin{table}
\parbox{.45\linewidth}{
\centering
\begin{tabular}{p{5mm}p{14mm}p{10mm}p{15mm}p{11mm}}
            \hline
            Id & Revenue & Served & Accepted & Missed \\ \hline
            
0 & 34.19 & 15.44 & 15.46 & 0.02 \\ \hline 
1 & 35.37 & 16.93 & 16.98 & 0.05 \\ \hline 
2 & 38.82 & 18.89 & 19.06 & 0.17 \\ \hline 
3 & 38.24 & 18.68 & 18.86 & 0.18 \\ \hline 
4 & 38.72 & 18.72 & 18.91 & 0.19 \\ \hline 
5 & 38.12 & 18.53 & 18.73 & 0.2 \\ \hline 
6 & 64.12 & 28.59 & 29.63 & 1.04 \\ \hline 
7 & 62.55 & 28.29 & 29.04 & 0.75 \\ \hline 
8 & 75.14 & 35.57 & 36.49 & 0.92 \\ \hline 
9 & 73.57 & 35.7 & 36.66 & 0.96 \\ \hline 
10 & 76.28 & 36.08 & 37.15 & 1.07 \\ \hline 
11 & 73.56 & 35.08 & 36.12 & 1.04 \\ \hline 
12 & 78.05 & 34.25 & 38.26 & 4.01 \\ \hline 
13 & 74.54 & 33.08 & 34.91 & 1.83 \\ \hline 
14 & 99.82 & 45.97 & 49.3 & 3.33 \\ \hline 
15 & 96.6 & 48.03 & 51.41 & 3.38 \\ \hline 
16 & 108.11 & 50.57 & 53.44 & 2.87 \\ \hline 
17 & 103.14 & 50.75 & 53.48 & 2.73
\\ \hline 

\end{tabular} \caption{Result of OPT+basis strategy on Gaussian distributed customers with Misspecified travel time distribution assumption.}
}
\hfill
\parbox{.45\linewidth}{
\centering
\begin{tabular}{p{5mm}p{14mm}p{10mm}p{15mm}p{11mm}}
            \hline
            Id & Revenue & Served & Accepted & Missed \\ \hline
            
0 & 28.12 & 13.04 & 14.9 & 1.86 \\ \hline 
1 & 26.15 & 12.1 & 12.91 & 0.81 \\ \hline 
2 & 36.69 & 17.57 & 18.23 & 0.66 \\ \hline 
3 & 35.36 & 17.77 & 18.45 & 0.68 \\ \hline 
4 & 37.29 & 17.59 & 17.83 & 0.24 \\ \hline 
5 & 37.02 & 18.11 & 18.37 & 0.26 \\ \hline 
6 & 42.64 & 16.34 & 20.5 & 4.16 \\ \hline 
7 & 38.64 & 13.98 & 15.31 & 1.33 \\ \hline 
8 & 57.82 & 23.33 & 26.46 & 3.13 \\ \hline 
9 & 54.57 & 25.59 & 28.94 & 3.35 \\ \hline 
10 & 66.29 & 30.63 & 34.3 & 3.67 \\ \hline 
11 & 64.18 & 29.25 & 31.25 & 2.0 \\ \hline 
12 & 50.01 & 17.97 & 22.18 & 4.21 \\ \hline 
13 & 44.57 & 19.34 & 22.46 & 3.12 \\ \hline 
14 & 67.61 & 24.97 & 29.05 & 4.08 \\ \hline 
15 & 64.02 & 31.94 & 39.19 & 7.25 \\ \hline 
16 & 79.78 & 32.67 & 37.95 & 5.28 \\ \hline 
17 & 77.84 & 35.09 & 40.01 & 4.92
\\ \hline 

\end{tabular} \caption{Result of OPT+basis strategy on Cluster distributed customers with Misspecified travel time distribution assumption.}
}
\end{table}

\begin{table}
\parbox{.45\linewidth}{
\centering
\begin{tabular}{p{5mm}p{14mm}p{10mm}p{15mm}p{11mm}}
            \hline
            Id & Revenue & Served & Accepted & Missed \\ \hline
            
0 & 25.35 & 10.84 & 12.4 & 1.56 \\ \hline 
1 & 22.76 & 10.62 & 11.55 & 0.93 \\ \hline 
2 & 35.26 & 16.68 & 17.79 & 1.11 \\ \hline 
3 & 33.95 & 15.65 & 16.27 & 0.62 \\ \hline 
4 & 37.34 & 17.45 & 17.85 & 0.4 \\ \hline 
5 & 36.86 & 18.21 & 18.71 & 0.5 \\ \hline 
6 & 38.69 & 13.72 & 16.25 & 2.53 \\ \hline 
7 & 35.37 & 12.85 & 14.11 & 1.26 \\ \hline 
8 & 54.57 & 21.6 & 24.58 & 2.98 \\ \hline 
9 & 51.69 & 23.39 & 26.59 & 3.2 \\ \hline 
10 & 63.53 & 27.6 & 30.5 & 2.9 \\ \hline 
11 & 61.06 & 25.76 & 27.33 & 1.57 \\ \hline 
12 & 45.05 & 17.0 & 22.84 & 5.84 \\ \hline 
13 & 41.49 & 16.11 & 18.41 & 2.3 \\ \hline 
14 & 62.98 & 24.66 & 31.17 & 6.51 \\ \hline 
15 & 62.61 & 27.4 & 32.42 & 5.02 \\ \hline 
16 & 77.3 & 29.63 & 33.53 & 3.9 \\ \hline 
17 & 73.97 & 32.82 & 37.71 & 4.89
\\ \hline 

\end{tabular} \caption{Result of OPT+basis strategy on Uniform distributed customers with Misspecified travel time distribution assumption.}
\label{last_full_table}
}
\end{table}
\end{APPENDICES}

\end{document}